\documentclass[12pt,a4paper]{article}
\usepackage{graphicx}
\usepackage{bm}
\usepackage{lscape}
\usepackage{amsmath}       

\usepackage{color}

\begin{document}

\title{Annotated square root computation in Liber Abaci and De Practica Geometrie by Fibonacci}

\author{Trond Steihaug\\ 
Department of Informatics\\ University of Bergen, Norway }
\maketitle

\abstract{
We study the square root computation by Leonardo Fibonacci (or Leonardo of Pisa) in his MSS Liber Abaci from c1202 and c1228 and De Practica Geometrie from c1220. We annotate a translation of Liber Abaci based on  transcripts from 1857 and 2020 and a translation from 2002 and a transcription of De Practica Geometrie from 1862 and a translation from 2008. We show that
Fibonacci is demonstrating the same method for all examples in the MSS and that this method deviates from the traditional description of the digit--by--digit method.  The description of the method used by Fibonacci is all verbal and summarized in tables for each square root example.
The manuscripts and transcription of the Latin texts are incomplete for some of the examples and the transcription and translation contains minor discrepancies and  some of the tables are incomplete and the missing digits are inserted.
}
\section{Introduction}

Leonardo Fibonacci or Leonardo of Pisa lived around 1170 to 1250. He was born in Italy, but educated in North Africa. He travelled widely in the Mediterranean until around 1200 when he settled in Pisa. He wrote several books of which copies of {\em Liber abaci} from around 1202 with a second edition from around 1228, {\em De Practica geometrie} from around 1220, {\em Flos} from around 1225, and {\em Liber quadratorum} from around 1225. The books Liber Abaci and De Practica Geometrie contains numerous examples computing the integer part of square root of natural numbers.
Fibonacci consider 16 examples, 5 in Chapter 14 of Liber Abaci and 11 in Chapter 2 of De Practica Geometrie:
The only example common in the two manuscripts is the square root of 12345. However, the technique
to compute the final digit is different in the two manuscripts. The only example where the technique is expanded compared is finding the square root of 7234. For this example Fibonacci do two modifications, first he multiplies 7243 with $10^4$ and finds the integer part of the square root of 7234 using the same technique as for the other examples in Liber Abacci and then evolve digit by digit of the two remaining digits in the integer part of $\sqrt{72430000}$.

Fibonacci is using a recursive technique where he assumes that only the last digit in the integer part of the square root needs to be determined. So to determine the integer part $\sqrt{12345}$ the technique  assumes that the integer part of the square root $n=123,\ a=\lfloor\sqrt{n}\rfloor$ and the residual $n-a^2$ are known. In his book on Liber abaci L\"{u}neburg \cite[p.258]{Luneburg1993} writes ``Hier wird man nun von Fibonacci wieder einmal \"{u}berrascht. Er entwickelt nicht Schritt f\"{u}r Schritt vor unseren Augen den Tertianeralgorithmus\footnote{Algorithm taught in High School mathematics} zur berechnung von $\lfloor\sqrt{N}\rfloor$, vielmehr sagt Fibonacci, man solle gem\"{a}ss den obigen Entwicklungen $\lfloor\sqrt{123}\rfloor$ berechnen.'' The recursive thinking is further evident in the example of $\lfloor\sqrt{9876543}\rfloor$ in
De Practica Geometrie where Fibonacci consider the root of a seven digit number but first finds the root of the first five digits. The method used by Fibonacci is not the Hindu method described by Datta and Singh \cite[p.169--175]{Datta1935} which is the method later used by Vi\'{e}te or the computational method of al-Nasaw\={\i} (c. 1011 – c. 1075) where $\sqrt{57342}$ is demonstrated \cite{Suter1906}. It is a common misconcept that Chapter 14 of Liber Abaci contains nothing significant not found in the Euclid's {\em Elements} \cite[p.10]{Sigler2002}\cite[p.71]{Devlin2011}. Further, to compute the integer part of the square root a further misconcept is that this is the Indian-Arabic algorithm \cite{Vogel1971}. Early contributions to the History of Mathematics focus mainly on computing the fractional part of the square root \cite[p.30]{Cantor1900}.

It is likely that Fibonacci has been well recognized during his lifetime for his pedagogical abilities which is stated in a decree from the town of Pisa in 1241 \cite{Caianiello2014}. It is shown that all examples computing the integer part of the square root are solved using the same algorithm suitable in an teaching environment.

The remaining part paper is organized as follows. In section \ref{sec:2} we introduce the notation and the technique used by Fibonacci to compute the integer part of the square root. The next two sections are the annotated examples in Liber Abaci and De Practica Geometrie. The last section \ref{sec:4} contains computing the fractional part of the square root from Liber Abaci.

\section{The integer part and notation}\label{sec:2}
The computation of the integer part of the square root by Fibonacci is a recursive procedure. If the
integer part of the square root has $k$ digits the Fibonacci assumes that the first $k-1$ digits are known
so the only digit to be determined will be the last. In the first part in this section, we derive several bounds and
approximations for the last digit.
For each example of computation of the integer part of the square root in the annotation, we reference the approximation or
bound being used.
\subsection{Integer part of the square root}
Let the decimal representation of a positive integer $N$ be
 \[ N=\beta_{2k-1}10^{2k-1} + \beta_{2(k-1)}10^{2(k-1)} + \cdots + \beta_1 10 + \beta_0, \]
 where \(\beta_{2k-1} + \beta_{2(k-1)}>0\) (not both equal 0). Since
 \[\left(\beta_{2k-1}10 + \beta_{2(k-1)}\right)10^{2(k-1)} \leq N < \left(\beta_{2k-1}10 + \beta_{2(k-1)}+1\right)10^{2(k-1)}, \]
 the integer part of the square root of $N$ has $k$ digits.
 Partition $N$ so that
 \(N= n10^2+\beta_1 10 + \beta_0,\quad n=\sum_{j=2}^{2k-1} \beta_j10^{j-2}. \) The integer part of $\sqrt{n}$ will have $k-1$ digits.
 Assume\footnote{$\lfloor\cdot\rfloor$ is the floor function. Let $x$ be a real number, then $\lfloor x\rfloor$ is the largest integer smaller or equal $x$.}  \[\lfloor\sqrt{n}\rfloor = \alpha_{k-1}10^{k-2} + \cdots + \alpha_2 10 + \alpha_1.\]
 Let $a=\lfloor\sqrt{n}\rfloor$ and consider
 \begin{equation*} \nonumber
   N-(10a+\alpha_0)^2 = (n-a^2)10^2+\beta_1 10 + \beta_0-2\cdot10a\alpha_0-\alpha_0^2  \\
 \end{equation*}
 Let $\alpha_0$ be largest integer in $\{0,1,\ldots,9\}$ so that
 \begin{equation}\label{eq:bound_1}
 (n-a^2)10^2+\beta_1 10 + \beta_0-2\cdot10a\alpha_0-\alpha_0^2 \geq 0
 \end{equation}
 then \(N-(10a+\alpha_0+1)^2<0\) and
 \begin{equation*}
(n-a^2)10^2+\beta_1 10 + \beta_0-2\cdot10a\alpha_0-\alpha_0^2- 2(10a+\alpha_0)  <1.
\end{equation*}
Since all terms on the left side of the inequality are integer, we have
\begin{equation}\label{eq:bound_2}
(n-a^2)10^2+\beta_1 10 + \beta_0-2\cdot10a\alpha_0-\alpha_0^2 \leq 2(10a+\alpha_0).
\end{equation}
The last digit $\alpha_0$ will be the smallest $\alpha_0$ so that (\ref{eq:bound_2}) holds.
A third bound on $\alpha_0$ is based on (\ref{eq:bound_1})
\begin{equation}\label{eq:bound_3}
\alpha_0\leq \frac{(n-a^2)10^2+\beta_1 10 + \beta_0}{2\cdot10a+\alpha_0}
\leq  \frac{(n-a^2)10^2+\beta_1 10 + \beta_0}{2\cdot10a}
\end{equation}
Approximations to this upper bound are
\begin{equation}\label{eq:bound_4}
\frac{(n-a^2)10+\beta_1}{2a},
\end{equation}
and
\begin{equation}\label{eq:bound_5}
\frac{n-a^2}{\lfloor\frac{2a}{10}\rfloor}.
\end{equation}
For the last example in De Practica Geometrie Fibonacci will be using even fewer decimal digits in the nominator and denominator of (\ref{eq:bound_5}) 
\(\frac{\lfloor(n-a^2)/10\rfloor}{\lfloor 2a/100\rfloor}.\)

\subsection{The annotation}

In Hughes' translation of De Practica Geometrie \cite[p.36]{Hughes2008} he demonstrates Fibonacci's square root computation of $\sqrt{864}$ in seven steps. Each of the first 6 steps corresponds to inserting one digit in a table.
\begin{description}
\item [1:] The root will be a two-digit number. The square of the first digit $\alpha_1$ is
immediately less than the first digit 8 ($=n$) of $N$. So place 2 $(=a=\alpha_1)$ under the 6.
\item [2:] Subtract $2^2$ from 8 and put the remainder 4 ($=n-a^2$) over the 8.
\item [3:] Double 2 and put it under 2.
\item [4:] An approximation to $\alpha_0$ is (\ref{eq:bound_4}), \(\alpha_0\approx \lfloor 46/4\rfloor =11\).
Largest single digit less than 11 is 9, probably the last digit $\alpha_0$ of the root.
Put 9 in the first place under 4.
\item [5:] Multiply 9 by 4 (under 2) and subtract 36 from 46 to get 10. Place 10 on
the diagonal.
\item [6:] Square 9, subtract 81 from 104 (on the diagonal), to get the remainder
23.
\item [7:] The root is verified by showing that $\alpha_0$ satisfies the second bound (\ref{eq:bound_2}), i.e.
\[ (n-a^2)10^2+\beta_1 10 + \beta_0-2\cdot10a\alpha_0-\alpha_0^2 = 23 \leq 2(10a+\alpha_0)=2\cdot29.\]

\end{description}

Hughes \cite[p.36]{Hughes2008} illustrates Fibonacci's computation with a sequence of tables 1:  to 6: in Figure \ref{fig:Hughes864}.
\begin{figure}[th!]
\begin{center}
\begin{tabular}{ccccccl}
1:&2:&3:&4:&5:&6:&\\
   $\begin{array}{ccc}  & & \\ & \\\bm 8&\bm 6&\bm4\\ & 2&\\ &&\\ \end{array}$
  &$\begin{array}{ccc}  & & \\4& \\\bm 8&\bm 6&\bm4\\ & 2&\\ &&\\ \end{array}$
  &$\begin{array}{ccc}  & & \\4& \\\bm 8&\bm 6&\bm4\\ & 2& \\ &4&\\ \end{array}$
  &$\begin{array}{ccc}  & & \\4& \\\bm 8&\bm 6&\bm4\\ & 2&9\\ &4&\\ \end{array}$
  &  $\begin{array}{ccc}1& & \\4&0\\\bm 8&\bm 6&\bm4\\ & 2&9\\ &4&\\ \end{array}$
  &  $\begin{array}{ccc}1& & \\4&0\\\bm 8&\bm 6&\bm4\\ & 2&9\\ &4&\\ \end{array}$& \\
  &&&&&&[23]
\end{tabular}
\end{center}
\caption{Notation used in \cite[p.36]{Hughes2008} to illustrate the square root of 864.} \label{fig:Hughes864}
\end{figure}

With the notation used in this paper, the corresponding evolution of tables is given in Figure \ref{fig:Ball} where the subscripts of the digits correspond to the numbered tables in Figures \ref{fig:Hughes864}. The subscripted digits were first used in \cite{Ball2010}.
The actual computation is shown in Figure \ref{fig:864} where the line numbers correspond to the subscripts Note that $4_2\bm6$ in line marked 4: in Figure \ref{fig:864} is the number 46 (in tens) where $4$ (in hundreds) is computed in line 2 and $\bm 6$ is the second digit in $864$. The sequence of tables in Figure \ref{fig:Hughes864} can readily be reconstructed from the table in Figure \ref{fig:Ball}.
\begin{figure}[th!]
\[\begin{array}{l}
0:\ (\beta_2,\beta_1,\beta_0)=(\bm 8,\bm 6,\bm 4)\\
1:\ \alpha_1^2 \leq  \beta_2 = \bm8,  \alpha_1=2\\
2:\ \beta_2-\alpha_1^2 = \bm 8-4 =4\text{ (in hundreds)}\\
3:\ 2\alpha_1 = 4\text{ (in tens)}\\
4:\ \alpha_0 \leq \lfloor\frac{4_2\bm 6}{4_3}\rfloor = 11, \alpha_0=9,\\
5:\ 4_2\bm6-2\alpha_1\alpha_0 = 10 \text{ (in tens)}\\
6:\ 1_5 0_5\cdot10+\beta_0-\alpha_0^2=104 - \alpha_0^2 =23\\
6:\ 10\cdot10+\beta_0-\alpha_0^2=23 \leq 2(10\alpha_1+\alpha_0)=2\cdot29
\end{array}\]
\caption{The actual computation of $\sqrt{864}$.} \label{fig:864}
\end{figure}

\begin{figure}[th!]
\begin{center}
 \begin{minipage}[t]{.25\linewidth}
\fbox{\begin{tabular}{ccl}
        &      &$\left(23_6\right.$   \\
  $1_5$ &      &         \\
  $4_2$ &$0_5$ &         \\
 {\bf 8}&{\bf 6}&{\bf 4} \\
        &$2_1$ & $9_4$   \\
        &$4_3$ &         \\
\end{tabular}}
\end{minipage}
\begin{minipage}[c]{.25\linewidth}
\includegraphics[height=3.5cm]{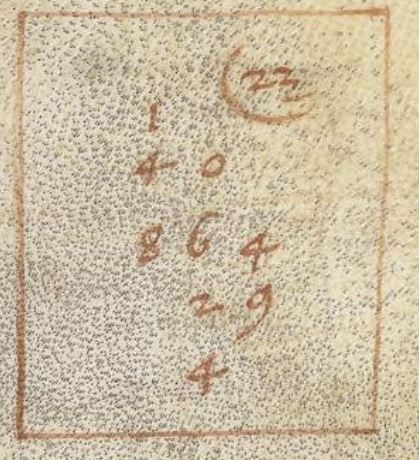}
\end{minipage}
\end{center}
\caption{Notation used in the paper and illustration from Folie 12r of MS  Urb. lat. 292 in  Città del Vaticano, Biblioteca Apostolica Vaticana.}\label{fig:Ball}
\end{figure}
 The table in Figure \ref{fig:Ball} is used by  Ball \cite{Ball2010}. The copy of Folie 12r of MSS  Città del Vaticano, Biblioteca Apostolica Vaticana, Urb. lat. 292 shows the actual table.\footnote{https://digi.vatlib.it/view/MSS\_Urb.lat.292. Accessed 01102023.}

The square root of 864 is 29 with remainder 23 which is marked as [23] in Figure \ref{fig:Hughes864} and $\left(23_6\right.$ in Figure \ref{fig:Ball}. Fibonacci's computation of the remainder is just an rearrangement of the binomial expression for $N-(\alpha_1 10+\alpha_0)^2$.

 \begin{alignat*}{2}
 N-(10\alpha_1 + \alpha_0)^2&=\beta_2 10^2+\beta_1 10+\beta_0-(10\alpha_1 + \alpha_0)^2&&\\
 &= \left(\left\{\left[(\beta_2 -\alpha_1^2) 10^2 +\beta_1 10\right]
 -(2\alpha_1)\alpha_0 10\right\}+\beta_0\right)-\alpha_0^2&&\\
 &=\left(\left\{\left[(4) 10^2 +\beta_1 10\right] -(2\alpha_1)\alpha_0 10\right\}+\beta_0\right)-\alpha_0^2&\quad& \text{by 2:}\\
 &=\left(\left\{\left[(4) 10 +\beta_1 \right]10 -(4)\alpha_0 10\right\}+\beta_0\right)-\alpha_0^2
 &\quad& \text{by 3:}\\
  &= \left(\left\{\left[46\right] -(4)\alpha_0 \right\}\cdot 10+\beta_0\right)- \alpha_0^2
   &\quad& \text{by 4:}\\
   &=\left(\left\{10 \right\} \cdot 10 +4\right)-\alpha_0^2    &\quad& \text{by 5:}\\
  &=(104)-9^2=23&\quad& \text{by 6:}.
 \end{alignat*}

 In the annotated transcription by Boncompagni \cite[p.19]{Boncompagni1862} of computing the integral part of the square root of 864, a reference [L:2] in the Latin text refers to line two (for each digit inserted in the table) in Figure \ref{fig:864} and Figure \ref{fig:Ball}, and [L:6.1] refers to the first of the lines marked 6: in Figure \ref{fig:864}.
 \begin{quote}
 Item si uis inuenire radicem de 864, pone 2 sub
6, cum 2 sint integra radix de 8 {\color{blue}[L:1, line 1 in Figure \ref{fig:864}]};
et 4 que remanent pone super 8 {\color{blue}[L:2]}:
deinde dupla ipsa
2, erunt 4, que pone sub 2 {\color{blue}[L:3]}; et per ipsa 4 diuide 46, scilicet copulationem superflui
terti\c{e} figure cum secunda, exibunt 11: ex qua diuisione possumus habere arbitrium
sequentis ponende figure, que multiplicanda est per duplum prime figure posite; et
postea per se ipsam, erit ipsa figura aut parum minus, aut totidem quantum ex ipsa
diuisione euenit; quod cognosces ex usu. Quare ponemus arbitrio 9 sub prima figura,
que sunt minus de 11 predictis  {\color{blue}[L:4]}; et multiplicabis 9 per 4, scilicet per duplum inuenti
binarij, et extrahes de 46 predictis, et ex remanentibus 10 pones 0 super 6, et 1
super 4 {\color{blue}[L:5]}; et copulabis ipsa 10 cum 4 primi gradus, erunt 104; de quibus extracta mul-
tiplicatione nouenarij in se ipso, remanent 23 {\color{blue}[L:6.1]}; que 23 sunt minus dupla radicis inuenct\c{e} {\color{blue}[L:6.2]}.
 \end{quote}

 A liberal translation based on the use of google translate\footnote{https://translate.google.com/?sl=la. Checked 1.11.2023.} and ChatGPT\footnote{https://chat.openai.com/auth/login. Checked 1.11.2023.} and the translation \cite[p.39]{Hughes2008} is given below. In the translation we have kept the Fibonacci's ordering of the digits which is numbering from the right, i.e. in the number 864 the digit 4 is the first, 6 is the second, and 8 is the third and last digit.

  \begin{quotation}{\em
 Suppose you wish to find the root of 864.}
 \begin{enumerate}
 \item {\em Put 2 under the 6 because 2 is the whole root of 8} {\color{blue}[L:1, line 1 in Figure \ref{fig:864}]}.
 \item {\em Put the remainder 4 above the 8} {\color{blue}[L:2]}.
 \item {\em Then double 2 to get 4 placing it under the 2 {\em{\color{blue}[L:3]}}.}
 \item {\em Form 46 from the 4 above the 8 and the 6 in the second place. Now divide the new number 46 by 4 to get 11. From this division we get an idea of the following first digit which must be multiplied by twice the digit you already found. Afterwards, square it. The digit is a little less or exactly as much as what comes from the division. Practice with this procedure will perfect you. So we choose 9 since it is less than 11, and put it under the first digit} {\color{blue}[L:4]}.
 \item {\em Multiply 9 by 4 which is two times the second digit and subtract the product from 46. The remainder is 10. Put 0 over the 6 and 1 above the 4 {\em{\color{blue}[L:5]}}}.
 \item {\em Join 10 with 4 in the first place to make 104. Subtract the square of 9 from it to get 23 {\em {\color{blue}[L:6.1]}} which is less than twice the root found} {\color{blue}[L:6.2]}.
     \end{enumerate}
 \end{quotation}

 More on notation in the paper:
 \begin{itemize}
  \item Boldface digits are the digits of the number to be taken square root of.
  \item The columns in the tables in Figure \ref{fig:Hughes864}  and \ref{fig:Ball} correspond to hundreds, tens, and ones.
  \item There are some changes in the notation in the revised {\em Liber Abaci} (from around 1228). Like in the computation of the integer part of $\sqrt{8754}$ in the revised edition 1 is written above 6 and 1 above 5 while in the first edition 11 is written above 5. We use the notation from the revised version here.
 \begin{figure}[th!]
 \centerline{\fbox{\begin{tabular}{cccl}
        &     &     &$\left(105_5\right.$   \\
        &$6_2$  &$11_4$ &         \\
 {\bf 8}&{\bf 7}&{\bf 5}&{\bf 4} \\
        &       &$9_1$ & $3_3$   \\
        &       &$9_1$ & $3_3$   \\
\end{tabular}}
      \fbox{\begin{tabular}{cccl}
        &       &       &$\left(105_5\right.$   \\
        &$1_4$  &       &         \\
        &$6_2$  &$1_4$  &         \\
 {\bf 8}&{\bf 7}&{\bf 5}&{\bf 4} \\
        &       &$9_1$ & $3_3$   \\
        &       &$9_1$ & $3_3$   \\
\end{tabular}}}
\caption{Minor notational differences in the two Liber Abaci editions}\label{fig:Notation}
\end{figure}
    \item Equal consecutive digits in a column in the MSS are compressed to one digit. With
         the subscripts introduced here, we have expanded the number of digits in a column with consecutive digits.
    \item In {\em De Practica Geometrie} a semicircle is placed in front of the remainder. In the paper we use ``$($" as shown in Figure \ref{fig:Ball}.
        The remainder (or residual) is not a part of the tables in Liber Abaci. However, in this paper we
        include it in upper right hand corner similar to in Practica Geometrie as in Figure \ref{fig:Notation}.
  \item In both editions of {\em Liber Abaci} the digits of the root are written two times
        (like $9_1$ and $3_3$ in Figure \ref{fig:Notation} to compute the sum), while in
        {\em De Practica Geometrie} the digits are multiplied by 2.
    \item Fibonacci calls  $\alpha_0$ the first digit.
    \item Red numbers in the tables are insertions due to errors, missing numbers, or for completeness.
    \item The blue inserted text in the transcribed and translated and  material are [L:$\cdot$]
    that refers to a line in the table for the square root computation or [C:$\cdot$] giving reference to additional information.
    \item A fraction like $3\frac{1}{6}$ is written  $\frac{1}{6}3$ by Fibonacci, but in this paper we use the standard notation $3\frac{1}{6}$.
    \end{itemize}

\section{Square root computation in Liber Abaci}\label{sec:LiberAbaci}
The examples of computing the integer part of the square root in Liber Abaci are
743, 8754, 12345, 927435, and 7234. The last examples deviates a little from the other four examples. The number is multiplied with $10^4$ and the integer part of the square root of $7234$ and the remainder is computed. The remaining two digits of the square root of 72340000 is explicitly stated and the remainder is computed.

\subsection{Liber Abaci $\sqrt{743}$ }
This is the first example of computing the integer part of the square root. The text contains two figures, the first showing computing the first digit in the root and its remainder putting the numbers in the correct columns corresponding to tens and hundreds.
Inserted text is bracketed and in color blue. A reference [L:2] refers to line two (for each number inserted in the table) in Figure \ref{fig:sqrt743} and [L:3.1] refers to the first of the lines marked 3:. This problem is discussed in \cite[p.142-144]{Pisano2015}. In this example $N=743$ and $n=7$.  To determine the last digit $\alpha_0$ inequalities (\ref{eq:bound_1}) and (\ref{eq:bound_2}) are used.

The following is the explanation on
how to complete the hand calculation given a more modern interpretation
The following explanation is based on the transcripts by Boncompagni \cite[p.353-354]{Boncompagni1857} and Giusti \cite[p.549--550]{Giusti2020}, the translation by Sigler \cite[p.491]{Sigler2002} and liberal use of google translate\footnote{https://translate.google.com/?sl=la. Checked 1.11.2023.} and ChatGPT\footnote{https://chat.openai.com/auth/login. Checked 1.11.2023.} is below.

\emph{
And if you wish to find the root of 743 using the {\em abaci} method, you can follow these steps:}
\begin{enumerate}
\item {\em Since 743 is a three-digit number, its root is a two-digit number. The last place of the root is taken below the second place, namely below the 4.
Twice put the largest root that the 7, namely the last digit of the 743, has in integers, below this 4.
The digit will be 2 which is twice put below the 4.} {\color{blue}[L:1, line 1 in Figure \ref{fig:sqrt743}]}
\item {\em Multiply one of the 2 by the other. There will be 4 which you subtract from the 7 leaving 3, and you put the 3 above the 7, as shown in Figure \ref{fig:sqrt743-1}} {\color{blue} [L:2]}.
\item {\em Couple the 3 with the preceding digit, namely with the 4, making 34. For this, you put twice some digit before the put ones, namely below the first place of the 743, and when it is multiplied by double the two and the product is subtracted from the above written 34, there remains a number which coupled with the first digit of 743, namely with the 3, can then be subtracted from the product of the digit put under the first place by itself} {\color{blue}[L:3.1, first line]}.
    {\em There remains then a number which does not exceed double the total of the found root} {\color{blue}[L:3.2]}, {\em and this will be the digit 7.
    The 7 is twice put below the 3} {\color{blue}[L:3.3]}.
\item {\em Multiply the upper 7 by the lower 2, and the lower 7 by the upper 2. You will thus have 28} {\color{blue}[L:4.1]}, {\em which subtracted from the 34 leaves 6. Put the 6 above the 4 of the 34} {\color{blue}[L:4.2]}.
    \item {\em Couple it with the preceding digit, namely the 3. There will be 63 from which you subtract the product of the 7 by the 7, namely 49. There remains 14.} {\color{blue}[L:5]}
        \end{enumerate}
{\em You will have 27 for the root of the 743; there remains 14, as is shown in the last illustration} {\color{blue}Figure \ref{fig:sqrt743}}.

\begin{figure}[th!]
\begin{center}
\fbox{\begin{tabular}{ccl}
        &      &\ \      \\
  $3_2$ &      &         \\
 {\bf 7}&{\bf 4}&{\bf 3} \\
        &$2_1$ &         \\
        &$2_1$ &         \\
\end{tabular}}
\caption{The first digit and the remainder in their correct column and row positions                                           in Fibonacci's example for  $\sqrt{743}$ in Liber Abaci}
 \label{fig:sqrt743-1}
 \end{center}
 \end{figure}
\begin{figure}[th!]
 \begin{center}

 \begin{minipage}[t]{.20\linewidth}
\fbox{\begin{tabular}{ccl}
        &      &$\left(14_5\right.$   \\
  $3_2$ &$6_4$ &         \\
 {\bf 7}&{\bf 4}&{\bf 3} \\
        &$2_1$ & $7_3$   \\
        &$2_1$ & $7_3$   \\
\end{tabular}}
\end{minipage}
\begin{minipage}[t]{.6\linewidth}
\(\begin{array}{l}
\text{1:}\ \alpha_1^2 \leq  \beta_2 = \bm 7,  \alpha_1=2\\
\text{2:}\ \beta_2-\alpha_1^2 = \bm 7-4=3 \text{ (in hundreds)}\\
\text{3:}\ (3_2{\bm 4}-2\alpha_1\alpha_0)\text{ (in tens))}, (34-4\alpha_0)\cdot 10+\beta_0\geq \alpha_0^2\\
\text{3:}\ ((34-4\alpha_0)\cdot 10+\beta_0- \alpha_0^2 \leq 2(10\alpha_1+\alpha_0),\\
\text{3:}\ \text{The two inequalities gives } \alpha_0=7\\
\text{4:}\ 2\alpha_1\alpha_0 =2\cdot2\cdot7=28\\
\text{4:}\ 34-2\alpha_1\alpha_0 = 6 \text{ (in tens)}\\
\text{5:}\ 6_4\cdot10 +\beta_0 - \alpha_0^2=63 - \alpha_0^2 =14
\end{array}\)
\end{minipage}
 \caption{The final figure in Fibonacci's example for  $\sqrt{743}$ in Liber Abaci.}
 \label{fig:sqrt743}
 \end{center}
 \end{figure}
Consider the inequalities in the unknown $\alpha_0$ in Figure \ref{fig:sqrt743}:
 \begin{eqnarray*}
 (34-4\alpha_0)\cdot 10+\beta_0- \alpha_0^2&=& 343-40\alpha_0-\alpha_0^2\geq 0, \text{ and}\\
 (34-4\alpha_0)\cdot 10+\beta_0- \alpha_0^2 - 2(10\alpha_1+\alpha_0)&=&303-42\alpha_0-\alpha_0^2 \leq 0.
 \end{eqnarray*}
 The above inequalities  are valid in the interval
 \[[2 \sqrt{186} - 21, \sqrt{743} - 20]\approx [6.27\ldots,7.25\ldots]\]
 and only integer solution is  $\alpha_0=7$. The two inequalities are property 2 and 3 in Pisano \cite[p.144]{Pisano2015}.  The square root computation is also found in \cite[p.157]{Friedlein1869}.

\subsection{
Liber Abaci $\sqrt{8754}$}

The following explanation is based on the transcripts by Boncompagni \cite[p.354]{Boncompagni1857} and Giusti \cite[p.550]{Giusti2020}, the translation by Sigler \cite[p.492]{Sigler2002} and use of google translate and ChatGPT--3.5.
A reference like [L:2] refers to line two (for each number inserted in the table) in Figure \ref{fig:sqrt8754} and [L:3.1] refers to the first of the lines marked three.  Here is $N=8754$ and $n=87$. The last digit is determined using (\ref{eq:bound_2}).

{\em If you wish to find the root of 8754 that is a number of four digits,
then we know that the root of it is a number of two digits;}
\begin{enumerate}
\item {\em Start with the number 87 and find the largest square number (a number that can be the result of multiplying itself) that is less than or equal to 87. In this case, it's 81, which is 9 squared, which you twice put below the 5} {\color{blue}[L:1, line 1 in Figure \ref{fig:sqrt8754}]}.

\item {\em Multiply 9 by 9, which is 81, and subtract that from 87, leaving you with 6 above the 7}
{\color{blue}[L:2]}.

\item {\em This coupled with
the preceding digit, namely the 5, makes 65 for which you put twice before the
put nines some figure, that when multiplied by double the 9, and the product
subtracted from the 65 leaves a number} {\color{blue}[L:3.1]}  {\em which when coupled with the digit of
the first place, namely the 4, can then be subtracted from the product of the
digit below the first place by itself and does not leave more than double the
entire found root, and the digit will be 3; this is twice put below the 4 before
the put nines} {\color{blue}[L:3.2]}.

\item {\em You will multiply crosswise the 3 by the 9 and the 3 by the 9;
there will be 54} {\color{blue}[L:4.1]} {\em that you subtract from the 65 leaving 11 that you put above the 65} {\color{blue}[L:4.2]}.

\item {\em Couple the 11 with the 4 that is in the first place; there will
be 114 from which you subtract the product of the 3 by the 3, namely 9; there
will remain 105; therefore the root of 8754 is in integers 93, and 105 remains} {\color{blue}[L:5]}.
\end{enumerate}

\begin{figure}[th!]
 \begin{center}

 \begin{minipage}[t]{.25\linewidth}
\fbox{\begin{tabular}{cccl}
        & $1_4$    &     &$\left(105_5\right.$   \\
        &$6_2$  &$1_4$ &         \\
 {\bf 8}&{\bf 7}&{\bf 5}&{\bf 4} \\
        &       &$9_1$ & $3_3$   \\
        &       &$9_1$ & $3_3$   \\
\end{tabular}}
\end{minipage}
\begin{minipage}[t]{.6\linewidth}
\(\begin{array}{l}
1:\ \alpha_1^2 \leq  10\beta_3+\beta_2 = \bm8\bm7,  \alpha_1=9\\
2:\ 10\beta_3+\beta_2-\alpha_1^2 = 87-9^2=6 \text{ (in hundreds)}\\
3:\ (6_2\bm5-2\alpha_1\alpha_0)\text{ (in tens)}, \\
3:\ (65-18\alpha_0)\cdot 10+\beta_0- \alpha_0^2\leq 2(10\alpha_1+\alpha_0), \alpha_0=3\\
4:\ 2\alpha_1\alpha_0= 2\cdot 3\cdot 9=54\\
4:\ 65-2\alpha_1\alpha_0 = 11 \text{ (in tens)}\\
5:\ 11_4\cdot10 +\beta_0 - \alpha_0^2=114 - \alpha_0^2 =105
\end{array}\)
\end{minipage}
 \end{center}
 \caption{Fibonacci's example for  $\sqrt{8754}$ from Liber Abaci}
 \label{fig:sqrt8754}

 \end{figure}
 Consider the inequality in line 3 in  Figure \ref{fig:sqrt8754}. This is inequality (\ref{eq:bound_2})
 \[(65-18\alpha_0)\cdot 10+\beta_0- \alpha_0^2- 2(10\alpha_1+\alpha_0)=474-182\alpha_0-\alpha_0^2 \leq 0\]
 and with $\alpha_0\geq 0$ is valid for \(\alpha_0\geq \sqrt{8755} - 91\approx 2.57.\) Smallest integer is $\alpha_0=3$.

 \subsection{ Liber Abaci $\sqrt{12345}$}
This is the only problem common for {\em Liber Abaci} and {\em De Practica Geomerie}. Here $N=12345$ and $n=123$.

The following translation is based on the transcripts by Boncompagni \cite[p.354--355]{Boncompagni1857} and Giusti \cite[p.550]{Giusti2020}, the translation by Sigler \cite[p.492]{Sigler2002} and use of google translate and ChatGPT--3.5.
A reference [L:2] refers to line two (for each number inserted in the table) in Figure \ref{fig:LA_sqrt12345} and [L:3.1] refers to the first of the lines marked 3:. The inequalty to determine the last digit is inequality (\ref{eq:bound_2}) while in De Practica Geomerie it is inequality (\ref{eq:bound_1}).

{\em If you wish to find the root of a number of five digits, say 12345,
you indeed find  the root of the number made of the
last three digits, namely the 123, and it will be 11, and 2 remains.}
\begin{enumerate}
\item {\em You therefore
twice put the 11 below the third and second places }{\color{blue}[L:1, line 1 in Figure \ref{fig:LA_sqrt12345}]}.
\item {\em There remains 2 that
you put above the 3} {\color{blue}[L:2]}.
\item {\em You couple it with the preceding digit, namely
the 4; there will be 24 which you put in the first place above the root, namely
before the put 11; you twice put a number that when multiplied crosswise by the
11, and the products are subtracted from the 24, leaves a number} {\color{blue}[L:3.1]} {\em which when
coupled with the digit of the first place, namely the 5, you subtract from it the
product of the digit and itself, and there does not remain more than double the
found root, and it will be 1 that is put before both of the 11} {\color{blue}[L:3.2]}.
\item {\em  You multiply
it crosswise by the 11; there will be 22 }{\color{blue}[L:4.1]} {\em that you subtract from the 24; there
remains 2 above the 4} {\color{blue}[L:4.2]}.
\item {\em This coupled with the 5 of the first place makes 25 from
which you subtract the product of the upper 1 and the lower 1; there remains
24} {\color{blue}[L:5]}, {\em and thus you will have a number of three digits, namely 111, for the root
of 12345, as should be, and there remains 24.}
\end{enumerate}

\begin{figure}[th!]
 \begin{minipage}[t]{.3\linewidth}
 \begin{center}
\fbox{\begin{tabular}{ccccl}
        &        &       &      &$\left(24_5\right.$  \\
        &        &$2_2$  &$2_4$ &        \\
 {\bf 1}&{\bf 2} &{\bf 3}&{\bf 4}&{\bf 5}\\
        &        &$1_1$ & $1_1$ & $1_3$  \\
        &        &$1_1$ & $1_1$ & $1_3$  \\
\end{tabular}}
\end{center}
\end{minipage}
\begin{minipage}[t]{.55\linewidth}
\(\begin{array}{l}
1:\   \lfloor \sqrt{\bm1\bm2\bm3}\rfloor = 11, \alpha_2=1, \alpha_1=1 \\
2:  \bm1\bm2\bm3-11_1^2=2  \text{ (in hundreds)}\\
3:\ 2_2{\bm 4}-2(\alpha_2\cdot10+\alpha_1)\alpha_0 \text{ (in tens)}, \\
3:\ (24-22\alpha_0)\cdot 10+\beta_0- \alpha_0^2\leq 2(11_1\cdot10+\alpha_0), \alpha_0=1\\
4:\ 2(\alpha_2\cdot10+\alpha_1)\alpha_0= 2\cdot 11\cdot 1=22\\
4:\ 24-2(\alpha_2\cdot10+\alpha_1)\alpha_0 = 2 \text{ (in tens)}\\
5:\ 2\cdot10 +\beta_0 - \alpha_0^2=2_4{\bm 5} - \alpha_0^2 =24
\end{array}\)
\end{minipage}
\caption{Fibonacci's example for  $\sqrt{12345}$ from Liber Abaci}
 \label{fig:LA_sqrt12345}

 \end{figure}

\subsection{
 Liber Abaci $\sqrt{927435}$ }
In this example $N=927435$, $n=9274$, $a=96$, and $n-a^2=58$. Inequality (\ref{eq:bound_2}) is used to determine
the last digit in the integer part of square root of $N$.

The following emphasized  text is based on the transcripts by Boncompagni \cite[p.355]{Boncompagni1857} and Giusti \cite[p.551]{Giusti2020}, the translation by Sigler \cite[p.492-493]{Sigler2002} and use of google translate and ChatGPT--3.5.
A reference [L:2] refers to line two (for each number inserted in the table) in Figure \ref{fig:sqrt927435} and [L:3.1] refers to the first of the lines marked 3:.
\begin{quote}
{\em If you wish to find the root of a six--digit number, as 927435,
which must have a root with three digits, then the last place digit is
put below the digit of the third place, namely below the 4; you find therefore
the root of the number made of the last four digits, namely the 9274, and this
you do according to that we demonstrated above in the finding of a root of a
number of four digits, and the root will be the number 96, and 58 remains.}
\begin{enumerate}
\item
{\em You therefore twice put the 96 below the third and second places }{\color{blue}[L:1, line 1 in Figure \ref{fig:sqrt927435}]}.
\item {\em You put
the 58 above the 74 of the 9274} {\color{blue}[L:2]},
\item {\em Then you couple the 58 with the preceding
digit, namely the 3 that is in the second place; there will be 583  for which
you put in the first place of the root, namely before the 96, twice such a digit
when multiplied crosswise by the 96} {\color{blue}[L:3.1]},{\em and the products subtracted from the 583
leave a number which when coupled with the first digit, namely the 5, and
subtracted from this is the figure multiplied by itself, there is not left a number
greater than double the found root, and it will be 3 that is put before both
of the 96 }{\color{blue}[L:3.2]}.
\item {\em  You multiply the 3 by the 96 crosswise; there will be 576 }{\color{blue}[L:4.1]}{\em  that you
subtract from the 583; there remains 7 above the 3} {\color{blue}[L:4.2]}.
\item {\em Which coupled with the 5
of the first place makes 75; from this you subtract the product of the 3 and the
3, namely 9; there remains 66 }{\color{blue}[L:5]}.
\end{enumerate}
\end{quote}

\begin{figure}[th!]
\begin{center}
\fbox{\begin{tabular}{cccccl}
        &       &       &       &        &$\left(66_5\right.$   \\
        &       &$5_2$  &$8_2$  & $7_4$        \\
 {\bf 9}&{\bf 2}&{\bf 7}&{\bf 4}&{\bf 3}&{\bf 5} \\
        &       &       &$9_1$ & $6_1$  &$3_3$  \\
        &       &       &$9_1$ & $6_1$  &$3_3$ \\
\end{tabular}}
\end{center}
\[\begin{array}{l}
1:\ \text{Insert the two digits of } \lfloor \sqrt{\bm{9274}}\rfloor =96 \text{ twice }\\
2:\ \text{Insert the two digits of } \bm{9274}-(9_16_1)^2=58\\
3:\ 5_28_2\cdot 10+\beta_1-2\cdot9_16_1\alpha_0 = 5_28_2\bm3-192\alpha_0 \text{ (in tens)}\\
3:\ (583 -192 \alpha_0)\cdot 10+\beta_0- \alpha_0^2\leq 2(9_16_1\cdot 10+\alpha_0), \alpha_0=3\\
4:\ 2\cdot 96\alpha_0= 2\cdot 96\cdot 3_3=576\\
4:\ 583-2\cdot96\alpha_0 =583-576= 7 \text{ (in tens)}\\
5:\ 7_4\cdot10 + \beta_0 - \alpha_0^2=7_4\bm5 - 9 =66
\end{array}\]
\caption{Fibonacci's example for  $\sqrt{927435}$ from Liber Abaci.}
\label{fig:sqrt927435}
 \end{figure}

\subsection{ Liber Abaci $\sqrt{7234}$ }
 Fibonacci computes $\sqrt{72340000}$ which is 100 times $\sqrt{7234}$. The technique of adding zeros to the radicand in order to obtain higher accuracy had already been done (see \cite[p.607]{Vogel1971}). The "another method" can also be that Fibonacci evolves several digits in the process and not only the last digit $\alpha_0$. Let
 \[72340000 = \sum_{j=0}^7\beta_j10^j,\quad
 \lfloor \sqrt{72340000}\rfloor = \sum_{j=0}^3 \alpha_j 10^j\]
 In this example Fibonacci first computes the integer part of the square root of $N=7234$, $\lfloor\sqrt{7234}\rfloor=85$ and the residual $ N-\lfloor\sqrt{7234}\rfloor^2=9$. To compute the final digit in $\lfloor\sqrt{7234}\rfloor$ inequality (\ref{eq:bound_2}) is used. This is illustrated in Figure \ref{fig:sqrt72340000_Giusti_1}. However, the individual digits in the columns are shifted compared to a table for the square root of $7234$. The next two digits in the square root $\lfloor \sqrt{72340000}\rfloor$ are just assumed to be $0$ and $5$. This is illustrated in Figure \ref{fig:sqrt72340000_Giusti_2}. The inserted numbers in the final table Figure \ref{fig:sqrt72340000_Giusti_3} are to compute the double of the root and the residual.

The following emphasized  text is based on the transcripts by Boncompagni \cite[p.355--356]{Boncompagni1857} and Giusti \cite[p.552-553]{Giusti2020}, the translation by Sigler \cite[p.492-493]{Sigler2002} and use of google translate and ChatGPT--3.5.
A reference [L:2] refers to line two (for each number inserted in the table) in figures Figure \ref{fig:sqrt72340000_Giusti_1} to Figure \ref{fig:sqrt72340000_Giusti_3} and [L:3.1] refers to the first of the lines marked 3:.

{\em
 I have 72340000; I shall show how to find the root by another
method; because, as it was said, the root of a number of a eight--digit number is a number
of four digits.}
\begin{enumerate}
\item {\em You therefore put the 8 below the 0 of the fourth place, as it is
the largest root smaller than the 72, namely the last two digits of the number
in integers,} {\color{blue}[L:1. Line 1 in Figure \ref{fig:sqrt72340000_Giusti_1}]}
\item {\em You multiply the 8 by itself; there will be 64; this subtracted
from the 72 leaves 8 which you put above the 2,} {\em {\color{blue}[L:2]}}
\item {\em and you know to couple them
and make 83, and you double the 8 put below the 0; there will be 16; the 6 of it
you put below the 8, and the 1 of it you put afterwards} {\color{blue}[L:3]}.
\item {\em
You find a digit that
multiplied by the 16 almost makes 83 but leaves a number which coupled with
the 4 of the following place, having subtracted from it the square of the digit
namely it multiplied by itself, and there is not left more than double the found
root, and the digit will be 5 that is put below the third place before the 8} {\color{blue}[L:4][C:1]}.
\item {\em You
multiply the 5 by the 1, namely by the last place of the 16; there will be 5 that
you subtract from the 8 that is above the 2; there remains 3 above the 8} {\color{blue}[L:5]}.
\item {\em
This
coupled with the following 3 makes 33 from which you take the multiplication
of the 5 by the 6; there will remain 3, namely that which is in the sixth
place }{\color{blue}[L:6][C:2]}.
\item {\em This coupled with the following 4 makes 34 from which you subtract the
square of the five, namely 25; there remains 9 above the 4} {\color{blue}[L:7]}.
\begin{figure}[th!]
\begin{center}
\fbox{\begin{tabular}{cccccccccl}
        & $3_5$ &       &       &       &        &      &       \\
        & $8_2$ & {\color{red}$3_6$}       & $9_7$      &       &        &      &      \\
 {\bf 7}&{\bf 2}&$\bm 3$&{\bf 4}&{\bf 0}&{\bf 0} &{\bf 0}&{\bf 0}\\
        &       &       &       &$8_1$  & $5_4$  &       &       \\
        &       &       &$1_3$  &$6_3$  &        &       \\
\end{tabular}}
\end{center}
\[\begin{array}{l}
1:\ \alpha_3^2 \leq \beta_7 10 + \beta_6 = 72, \alpha_3=8\text{ (in thousands)}\\
2:\  \beta_7 10 + \beta_6 - \alpha_3^2 = 72-64=8\text{ (in millions) }\\
3:\ 8\cdot 10 + \beta_5=8\cdot 10 + \bm 3 =83, 2\alpha_3=16 \text{ (in thousands)}\\
4:\ (8\cdot 10 + \beta_5)10 + \beta_4 - 2\alpha_3\cdot 10\alpha_2 - \alpha_2^2 \leq 2(80+\alpha_2)\\
5:\ 8_2-5_4\cdot 1_3=3\\
6:\ 8_2\bm 3-2\alpha_3\alpha_2 = 83-16\cdot 5=(8_2-5_4\cdot1_3)\cdot10+\bm3-5_4\cdot6_3=3\\
7:\ 3_6\bm 4-\alpha_2^2=9
\end{array}\]
\caption{First illustration to compute the first two digits of $\sqrt{72340000}$ from Liber Abaci (2nd ed) \cite[p.552]{Giusti2020}.}
\label{fig:sqrt72340000_Giusti_1}
 \end{figure}
\item {\em You double
the 5; there will be 10; you put of it the 0 below the 5, and the 1 you add
to the 6 that is below the 8, and thus you will have 170 for the double of the
85 }{\color{blue}[L:8. Line marked L:8 in Figure \ref{fig:sqrt72340000_Giusti_2}]}.
\item {\em Next the 0 is put below the 0 of the second place before the put 85.}
{\color{blue}[L:9] [C:3]}
\item {\em The
multiplication of the second place by the fifth, namely by the 1, makes the sixth
place; this place is not there as the last digit of the remaining number, namely
the 9, is in the fifth place; this 0 put, you double it making 0 which you put
below it, namely before the 170, and you will have 1700 for double the found
root by which the put digit is multiplied} {\color{blue}[L:10]}.
\item {\em Therefore you put the 5 below the
0 in the first place {\em {\color{blue}[L:11][C:4]}}.}
\item {\em
You multiply it by the 1, and you subtract the product
from the 9; there remains 4 above it {\em {\color{blue}[L:12]}}.}
\item {\em You multiply the 5 by the 7, and you
subtract the product from the 40; there remains 5 above the 0 in the fourth
place {\em {\color{blue}[L:13]}}}.
 \begin{figure}[th!]
\begin{center}
\fbox{\begin{tabular}{cccccccccl}
        & $3_5$ &       &$4_{12}$  &       &        &      &       \\
        & $8_2$ &{\color{red}$3_6$}       &$9_7$  & $5_{13}$ &        &      &   \\
 {\bf 7}&{\bf 2}&$\bm 3$&{\bf 4}&{\bf 0}&{\bf 0} &{\bf 0}&{\bf 0}\\
        &       &       &       &$8_1$  & $5_4$  &$0_9$ & $5_{11}$ \\
        &       &       &$1_8$  &$7_8$  &$0_8$   &$0_{10}$&       \\
\end{tabular}}
\end{center}
\[\begin{array}{l}
8:\ 2\alpha_2=2\cdot5_4, 1_36_3+1_8= 1_87_8. \quad 1_87_80_8=2(\alpha_310+\alpha_2)=2\cdot85\\
9:\  \alpha_1= 0_9\\
10:\ 0_{10}=2\alpha_1 = 2\cdot0_9.\quad
1_6 7_6 0_6 0_{10}=2(\alpha_310^2+\alpha_2\cdot 10 + \alpha_1)=2\cdot850\\
11:\  \alpha_0= 5_{11}\\
12:\ 9_7-5_{11}\cdot 1_{8}= 4_{12}\\
13:\ 4_{12}\bm 0-5_{11}\cdot7_{10}=5_{13}
\end{array}\]
\caption{Second illustration of $\sqrt{72340000}$ from Liber Abaci}
\label{fig:sqrt72340000_Giusti_2}
 \end{figure}
\item {\em You multiply the 5 by the 0 which is below the 5, and you subtract
the product from the 50; there remains 50, and you multiply the 5 by the 0
which is below the 0 in the number 1700, and you subtract from the 500; there
remains 500 ending above the second place, and you multiply the 5 by itself,
and you subtract the product from the 5000; there remains 4975 above the 5000
{\em {\color{blue}[L:14. Line marked L:14 in Figure \ref{fig:sqrt72340000_Giusti_3}]}}},
\item {\em You double the 5, namely that which is in the first place of the found root;
there will be 10; from the 10 you put the 0 below the 5, and you put the 1 after
it, deleting the 0 which is in the place, and thus you will have 17010 for double
the found root {\em {\color{blue}[L:15]}},
as clearly shown in the third illustration {\em {\color{blue}[See Figure \ref{fig:sqrt72340000_Giusti_3}]}}.}
\end{enumerate}
{\em The root is 8505 and there remains 4975.}
 \begin{figure}[th!]
\begin{center}
\fbox{\begin{tabular}{cccccccccl}
        & $3_5$ &       &$4_{12}$&$4_{14}$&        &        &        \\
        & $8_2$ &{\color{red}$3_6$}       &$9_7$   &$5_{13}$&$9_{14}$&$7_{14}$&$5_{14}$\\
 {\bf 7}&{\bf 2}&$\bm 3$&{\bf 4} &{\bf 0}&{\bf 0}&{\bf 0}&{\bf 0}\\
        &       &       &        &$8_1$   &$5_4$   &$0_9$   &$5_{11}$ \\
        &       &       &$1_{15}$&$7_{15}$&$0_{15}$&$1_{15}$&$0_{15}$       \\
\end{tabular}}
\end{center}
\[\begin{array}{l}
14:\ 5_{13}\bm0 - 5_{11}\cdot0_{8}=50,\ \ 50\bm0-5_{11}\cdot0_{10}=500,\ \
500\bm0-5_{11}^2=4_{14}9_{14}7_{14}5_{14}\\
15:\ 1_6 7_6 0_6 1_{10} 0_{10} =2(\alpha_310^3+\alpha_2\cdot 10^2 + \alpha_1 10+\alpha_0)=2\cdot8505
\end{array}\]
\caption{Third illustration}
\label{fig:sqrt72340000_Giusti_3}
 \end{figure}

The most common deviation between the figures in the different manuscripts, transcripts and translations are in misalignment in the columns.  In Figure \ref{fig:sqrt72340000} the different figures from the square root computation in two manuscripts \cite{Fibonacci_LA_XI21}\footnote{https://archive.org/details/conventi-soppressi-c.-i.-2616. Accessed 01102023.} and \cite{Fibonacci_LA_2626}\footnote{https://archive.org/details/magliabechiano-xi.-21. Accessed 01102023.}, two transcripts \cite{Giusti2020,Boncompagni1857} and  the translation \cite{Sigler2002}.

 \begin{figure}[th!]
\begin{center}
\fbox{\begin{tabular}{cccccccl}
        &       &       &        &        &        &        &$\left(4975_{14}\right.$ \\
        & $3_5$ &$   $  &$4_{12}$&$4_{14}$&        &        &        \\
        & $8_2$ &{\color{red}$3_6$}       &$9_7$   &$5_{13}$&$9_{14}$&$7_{14}$&$5_{14}$\\
 {\bf 7}&{\bf 2}&$\bm 3$&{\bf 4} &{\bf 0} &{\bf 0} &{\bf 0} &{\bf 0}\\
        &       &       &        &$8_1$   & $5_4$  &$0_9$   &$5_{11}$ \\
        &       &       &$1_3$   &$6_3$   &        &        &       \\
        &       &       &$1_8$   &$7_8$   &  $0_8$ &$0_{10}$&       \\
        &       &       &$1_{15}$&$7_{15}$&$0_{15}$&$1_{15}$&$0_{15}$       \\
\end{tabular}}
\end{center}
\caption{A complete illustration of $\sqrt{72340000}$ based on Liber Abaci}
\label{fig:sqrt72340000_1}
 \end{figure}

Also note that in this example Fibonacci the last line in the table is two times the partial root. Further, a line in the table is erased/overwritten.
\begin{description}
\item[C:1] The inequality (\ref{eq:bound_2}):
\(674-162\alpha_2 - \alpha_2^2 \leq 0, \text{ smallest }\alpha_2=5. \)
\item[C:2] There is already a $\bm 3$ in the sixth position (counting from right) so this is reused for hand calculation. A new $3$ is inserted and marked with a subscript $3_6$ for clarity.
\item[C:3] Fibonacci gives no indication on how $\alpha_1$ is computed, but the residual $7234-85^2=9_6$ and the next two digits in $N$ are $\beta_3=\beta_2=\bm 0$. For $n=7234$, $a=85$ in (\ref{eq:bound_1}) the largest integer $\alpha_1$ will satisfy
\[9\cdot10^2-2\cdot85\cdot10\cdot\alpha_1-\alpha_1^2\geq 0 \text{ so } \alpha_1= 0=0_9.\]
\item[C:4] The residual $723400-850^2=9_700$ and the next two digits in $N$ are $\beta_1=\beta_0=\bm 0$,  For $n=723400$, $a=850$ in (\ref{eq:bound_1}) the largest integer $\alpha_0$ will satisfy
\[900\cdot10^2-2\cdot850\cdot10\cdot\alpha_0-\alpha_0^2\geq 0 \text{ so } \alpha_0= 5=5_{11}\]

\end{description}
\begin{figure}[p]
 \begin{center}
 \begin{minipage}[t]{.3\linewidth}
\includegraphics[scale=.3]{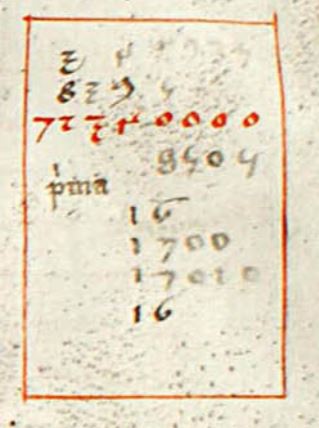}
\includegraphics[scale=.6]{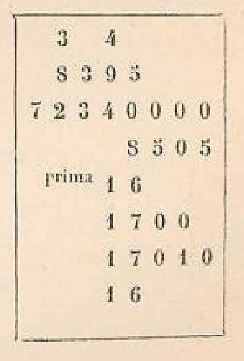}
\includegraphics[scale=.8]{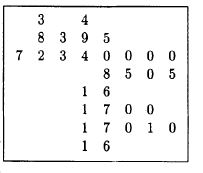}
\includegraphics[scale=.6]{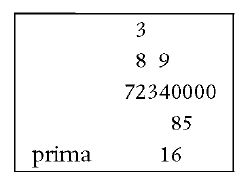}
\includegraphics[scale=.6]{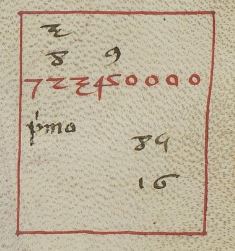}
\end{minipage}
\begin{minipage}[t]{.3\linewidth}
\includegraphics[height=3cm]{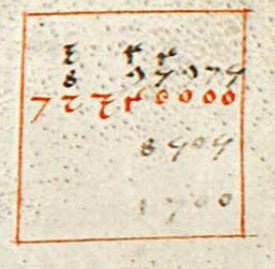}
\includegraphics[scale=.6]{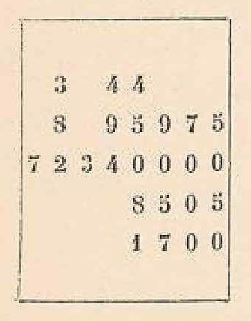}
\includegraphics[scale=.8]{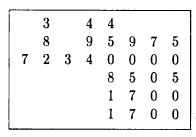}
\includegraphics[scale=.8]{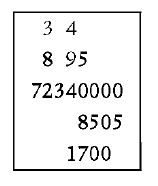}
\includegraphics[scale=.6]{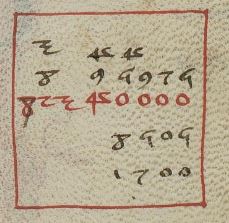}
\end{minipage}
\begin{minipage}[t]{.3\linewidth}
\includegraphics[height=3cm]{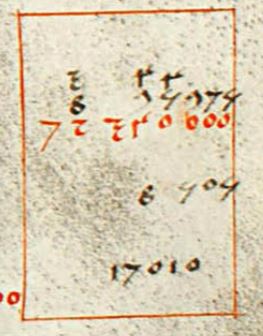}
\includegraphics[scale=.6]{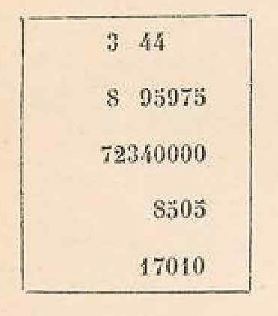}
\includegraphics[scale=.8]{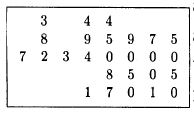}
\includegraphics[scale=.8]{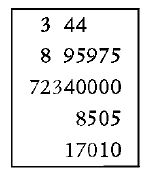}
\includegraphics[scale=.6]{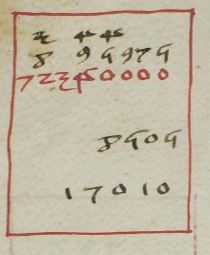}
\end{minipage}

 \caption{The three illustrations of $\sqrt{72340000}$ from five different sources. First row from BNCF, Conventi Soppressi C.I.2616, the second row from Boncompagni transcript from 1857, the third row from Sigler 2002, the fourth row from Giusti \cite[p.552]{Giusti2020} and the final row from BNCF, Magliabechiano XI.21  }
 \label{fig:sqrt72340000}
 \end{center}
 \end{figure}

\section{Square root computation in De Practica Geometrie}
De Practica Geometrie contains 11 examples of computing the integer part of the square root of the numbers 153, 864, 960, 1234, 6142, 8172,12345, 98765, 123456, 987654 and 9876543. Fibonacci makes some minor changes the notation compared to Liber Abaci. The text is based on the transcripts of Boncompagni from 1862 \cite{Boncompagni1862} and the translation by Hughes from 2008 \cite{Hughes2008}.
\subsection{
De Practica Geometrie $\sqrt{153}$ }
The emphasized text is based on the transcript of Boncompagni \cite[p.19]{Boncompagni1862} and the translation \cite[p.39]{Hughes2008}. In the computation of the square root of 153 and 854, the value of the last digit $\alpha_0$ is stated explicitly and at the end shown to satisfy (\ref{eq:bound_1}).
\begin{quote}{\em
For example, we wish to find the root of 153.}
\begin{enumerate}
\item {\em  You will find the root of 1 in
the third place to be 1. Place it under the 5 {\em{\color{blue}[L:1, line 1 in Figure \ref{fig:sqrt153}]}}}
\item {\em and put 2 before it under the 3 {\em{\color{blue}[L:2]}}.}
\item {\em  Multiply 2 by 2, which is twice the root that you found {\em{\color{blue}[L:3.1]}}, to get 4 {\em{\color{blue}[L:3.2]}}.}
    \item {\em Subtract this from 5 and put the remainder 1 above the 5 {\em{\color{blue}[L:4]}}.}
    \item {\em Join the 1 to the 3
in the first place to form 13 {\em{\color{blue}[L:5.1]}}. Subtract from this the square of 2 to get 9 {\em{\color{blue}[L:5.2]}} that
is less than 24 {\em{\color{blue}[L:5.3]}}, twice the root which has been found {\em{\color{blue}[L:5.4]}}.}
\end{enumerate}
{\em Consequently, the whole
root of 153 is 12 with remainder 9.}
\end{quote}
\begin{figure}[ht!]
\begin{minipage}[t]{.25\linewidth}
\fbox{\begin{tabular}{ccl}
        &      &$\left(9_5\right.$   \\
        &$1_4$ &         \\
 {\bf 1}&{\bf 5}&{\bf 3} \\
        &$1_1$ & $2_2$   \\
        &$2_3$ &          \\
\end{tabular}}
\end{minipage}
\begin{minipage}[t]{.6\linewidth}
\(\begin{array}{l}
1:\ \alpha_1^2 \leq  \beta_2 = {\bm 1},  \alpha_1=1\\
2:\ \alpha_0=2\\
3:\ 2\alpha_1=2\\
3:\ 2\alpha_1\alpha_0 = 2_3\alpha_0 =2_3\cdot 2_2 = 4\\
4:\ \beta_2-\alpha_1^2 = 1-1=0 \text{ (in hundreds)}\\
4:\ 0\cdot 10 +\beta_1-2\alpha_1\alpha_0 =\bm5-4=1 (\text{ in tens}) \\
5:\ (\bm5-2\alpha_1\alpha_0)\cdot 10+\beta_0=10+\bm3=13\\
5:\ (\bm5-2\alpha_1\alpha_0)\cdot 10+\beta_0 -\alpha_0^2=13-4=9\\
5:\ 2(10\alpha_1+\alpha_0)=24\\
5:\ (5-2\alpha_1\alpha_0)\cdot 10+\beta_0- \alpha_0^2=9\leq 2(10\alpha_1+\alpha_0)=24.

\end{array}\)
\end{minipage}
 \caption{Fibonacci's example for  $\sqrt{153}$ from De Practica Geometrie.}
 \label{fig:sqrt153}
 \end{figure}
\subsection{De Practica Geometrie $\sqrt{864}$ }
The emphasized text below is based on the transcript of Boncompagni \cite[p.19]{Boncompagni1862} and the translation \cite[p.39]{Hughes2008}.
The inserted text in the quote refers to  Figure \ref{fig:sqrt864}.
 \begin{quotation}{\em
 If you wish to find the root of 864,}
 \begin{enumerate}
 \item {\em put 2 under the 6 because 2 is the whole root of 8} {\color{blue}[L:1, line 1 in Figure \ref{fig:sqrt864}]}.
 \item {\em Put the remainder 4 above the 8} {\color{blue}[L:2]}.
 \item {\em Then double 2 to get 4 placing it under the 2 {\em{\color{blue}[L:3]}}.}
 \item {\em Form 46 from the 4 above the 8 and the 6 in the second place. Now divide the new number 46 by 4 to get 11. From this division we get an idea of the following first digit which must be multiplied by twice the digit you already found. Afterwards, square it. The digit is a little less or exactly as much as what comes from the division. Practice with this procedure will perfect you. So we choose 9 since it is less than 11, and put it under the first digit} {\color{blue}[L:4]}.
 \item {\em Multiply 9 by 4 which is two times the second digit and subtract the product from 46. The remainder is 10. Put 0 over the 6 and 1 above the 4 {\em{\color{blue}[L:5]}}.
 \item Join 10 with 4 in the first place to make 104. Subtract the square of 9 from it to get 23 {\em {\color{blue}[L:6.1]}} which is less than twice the root found} {\color{blue}[L:6.2]}\footnote{In \cite[p.40]{Hughes2008} this is stated as {\em This is less than 29 the root that has been found}}.
     \end{enumerate}
 \end{quotation}
 In Figure \ref{fig:sqrt864} consider
$\sqrt{N}$, where $N$ is the form  $N= \beta_2 10^2+\beta_1 10+\beta_0$ where  $N=864$. Let the integral part of $\sqrt{N}$ be $\alpha_1 10 + \alpha_0$.
 \begin{figure}[ht!]
 \begin{center}

 \begin{minipage}[t]{.2\linewidth}
\fbox{\begin{tabular}{ccl}
        &      &$\left(23_6\right.$   \\
  $1_5$ &      &         \\
  $4_2$ &$0_5$ &         \\
 {\bf 8}&{\bf 6}&{\bf 4} \\
        &$2_1$ & $9_4$   \\
        &$4_3$ &         \\
\end{tabular}}
\end{minipage}
\begin{minipage}[t]{.5\linewidth}
\(\begin{array}{l}
0:\ (\beta_2,\beta_1,\beta_0)=(\bm 8,\bm 6,\bm 4)\\
1:\ \alpha_1^2 \leq  \beta_2 = \bm8,  \alpha_1=2\\
2:\ \beta_2-\alpha_1^2 = \bm 8-4 =4\text{ (in hundreds)}\\
3:\ 2\alpha_1 = 4\text{ (in tens)}\\
4:\ \alpha_0 \leq \lfloor\frac{4_2\bm 6}{4_3}\rfloor = 11, \alpha_0=9,\\
5:\ 4_2\bm6-(2\alpha_1)\alpha_0 = 4_2\bm6-4_39_4 =10 \text{ (in tens)}\\
6:\ 1_5 0_5\cdot10+\beta_0-\alpha_0^2=104 - \alpha_0^2 =23\\
6:\ 10\cdot10+\beta_0-\alpha_0^2=23 \leq 2(10\alpha_1+\alpha_0)=2\cdot29
\end{array}\)
\end{minipage}
 \caption{Fibonacci's example for  $\sqrt{864}$ in De Practica Geometrie.}
 \label{fig:sqrt864}
 \end{center}
 \end{figure}

In Figure \ref{fig:sqrt864} the $4_2$ and {\bf 6} from the diagonal forms $46$ and is $\beta_2 10^2+ \beta_1 10 -(10\alpha_1)^2=100(\beta_2-\alpha_1^2)+\beta_1 10=460$ (or 46 since the digit 6 is placed in column 2). The next diagonal is $1_5$, $0_5$, and {\bf 4} forming 104. In Figure \ref{fig:sqrt864}  follows the detailed derivation in  \cite[p.36]{Hughes2008}.

Consider the inequality  in [L:6.2] in Figure \ref{fig:sqrt864}
\begin{eqnarray*}
 (46-4\alpha_0)\cdot 10+\beta_0- \alpha_0^2 - 2(10\cdot2+\alpha_0)&=&424-42\alpha_0-\alpha_0^2 \leq 0
 \end{eqnarray*}
 which for nonnegative $\alpha_0$ is valid for $\alpha_0\geq\sqrt{865} - 21 \approx 8.41\ldots$. Smallest integer is $\alpha_0=9$.

\subsection{De Practica Geometrie $\sqrt{960}$}
The emphasized text is based on the transcript of Boncompagni \cite[p.19--20]{Boncompagni1862} and the translation \cite[p.40]{Hughes2008}. The references in the text are to Figure \ref{fig:sqrt960}.
 \begin{quote}{\em
 Again, if you wish to find the root of 960,}
 \begin{enumerate}
 \item {\em put 3 the root of 9 under the 6 {\em{\color{blue}[L:1]}}}
 \item {\color{blue}[L:2]}\footnote{Added $0_2$ above {\bf 9} to be consistent.}
 \item {\em Double 3 to get 6 and place it under the 3 {\em{\color{blue}[L:3]}}}
 \item {\em  Now multiply 6 by the first digit  {\em{\color{blue}[L:4.1]}}.
Subtract the product from the upper 6 from which a digit remains {\em{\color{blue}[L:4.2]}}.
Join it with 0 {\em{\color{blue}[L:4.3]}}.
You can subtract the square of the digit {\em{\color{blue}[L:4.4]}}.
There will
be no remainder beyond twice the root that had been found. And that digit
is 0 {\em{\color{blue}[L:4.5]}}.}
\item {\em  Now the product of 0 with 6 (twice 3) subtracted from 6 leaves the
same 6 {\em{\color{blue}[L:5.1]}},
Join this with 0 in the first place to make 60. Subtracting
the square of 0 from this leaves 60, twice 30 the root that was found {\em{\color{blue}[L:5.2]}}.}
\end{enumerate}
\end{quote}

\begin{figure}[ht!]
 \begin{minipage}[t]{.2\linewidth}
\begin{center}
\fbox{\begin{tabular}{ccl}
  ${\color{red}0_2}$      &       &$\left(60_5\right.$   \\
 {\bf 9}&{\bf 6}&{\bf 0} \\
        &$3_1$  & $0_4$   \\
        &$6_3$  &       \\
\end{tabular}}
\end{center}
\end{minipage}
\begin{minipage}[t]{.5\linewidth}
\(\begin{array}{l}
1:\ \alpha_1^2\leq \bm 9, \alpha_1=3_1\\
2:\ \beta_2-\alpha_1^2=0_2\\
3:\ 2\alpha_1 = 6_3 \text{ (in tens)}\\
4:\ 2\alpha_1 \alpha_0=6_3\alpha_0\\
4:\ \bm 6- 6_3\alpha_0 \text{ (in tens) }\\
4:\ (6-6\alpha_0)10+\beta_0=60-60\alpha_0\\
4:\ (6-6\alpha_0)10+\beta_0-\alpha_0^2=60-60\alpha_0-\alpha_0^2\\
4:\ 60-60\alpha_0-\alpha_0^2 \leq 2(10\alpha_1+\alpha_0)\\
\quad\quad\quad=2\cdot(3_1\cdot10+\alpha_0)=60+2\alpha_0, \alpha_0=0_4\\
5:\ 0_2\cdot 10+\beta_1-2\alpha_1\alpha_0 = \bm6-6_3\cdot0_4=6\text{ (in tens) }\\
5:\ (\beta_1-2\alpha_1\alpha_0)\cdot 10 + \beta_0 = 60=2(10\alpha_1+\alpha_0)=2\cdot30

\end{array}\)
\end{minipage}
 \caption{Fibonacci's example for  $\sqrt{960}$ from De Practica Geometrie.}
 \label{fig:sqrt960}
 \end{figure}
\subsection{
De Practica Geometrie $\sqrt{1234}$ }
The emphasized text is based on the transcript of Boncompagni \cite[p.20]{Boncompagni1862} and the translation by Hughes \cite[p.40-41]{Hughes2008}. See Figure \ref{fig:sqrt1234} for the corresponding line numbers. This is also the first example where the detailed description of a part of the computation.
\begin{quotation}{\em
If you wish to know the root of a four digit number, find first the root
of the last two digits. Join the remainder with the remaining two digits
and proceed according to what we said above for three digit numbers. For
instance, we want to find the root of 1234. }
\begin{enumerate}
\item {\em Put 3 the root
of the square less than 12 under the 3 {\em{\color{blue}[L:1]}}}
\item {\em  and the remainder 3 over the 2 of 12
 {\em{\color{blue}[L:2]}}.}
\item {\em
Join 3 with the following digits to make 334 {\em{\color{blue}[L:3]}}.}
\item {\em  Then double the root which
you found namely 6 and place it under the root 3 {\em{\color{blue}[L:4]}}.}
\item {\em Now think: how many
multiplications are to be made according to the instructions given above
and how many digits will there be in the number from which we must subtract
those multiplications? There will be two multiplications. The first is
the product of the first digit by twice the [partial] root or 6. The second
is the square of the first digit {\em{\color{blue}[C:1]}}. Both of these products are to be subtracted
from 33 in order to finish the last multiplication under the digit in the first
place. Hence, the first product will be subtracted from 33 the union of the
last two digits, and the other product from the union of the remainder and
the 4 under the first digit {\em{\color{blue}[C:2]}}. Consequently, put 5 under the first digit because
33 divided by 6 leaves 5 {\em{\color{blue}[L:5]}}. }
\item {\em
Multiply 5 by 6 and subtract the product from 33
for a remainder of 3 in the second place {\em{\color{blue}[L:6]}}.}
\item {\em Join 3 with 4 in the first place to
make 34. From this subtract the square of 5 to leave 9 {\em{\color{blue}[L:7]}}}.
\end{enumerate}
{\em Thus, you have 35
as the root of 1234 with remainder 9.}
\end{quotation}
\begin{figure}[ht!]
 \begin{minipage}[t]{.3\linewidth}
\begin{center}
\fbox{\begin{tabular}{cccl} &       &       &$\left(9_7\right.$   \\
        &$3_2$  &{\color{red}$3_6$}       &        \\
 {\bf 1}&{\bf 2}&{\bf 3}&{\bf 4} \\
        &       &$3_1$  & $5_5$   \\
        &       &$6_4$  &       \\
\end{tabular}}
\end{center}
\end{minipage}
\begin{minipage}[t]{.5\linewidth}
\(\begin{array}{l}
1:\ \alpha_1^2 \leq  \beta_3 10+\beta_2= \bm{12},  \alpha_1=3\\
2:\ \beta_3 10+\beta_2-\alpha_1^2=\bm{12}-9=3 \text{ (in hundreds)}\\
3:\ 3_2\cdot100 +\beta_1\cdot10 + \beta_0=3_2\bm{34}\\
4:\ 2\alpha_1 =6 \text{ (in tens)}\\
5:\ \alpha_0 \approx \lfloor\frac{3_2\bm3}{6}\rfloor = 5, \alpha_0=5,\\
6:\ 3_2{\bm 3}-2\alpha_1\alpha_0 = 3 \text{ (in tens)}\\
7:\ 3_6\cdot10+\beta_0-\alpha_0^2=3_64 - \alpha_0^2 =9
\end{array}\)
\end{minipage}
 \caption{Fibonacci's example for  $\sqrt{1234}$ from De Practica Geometrie}
 \label{fig:sqrt1234}
 \end{figure}
 There is a $3$ in the second place in Figure \ref{fig:sqrt1234}, so the $\color{red}3_6$ is added for clarity to the $\bm 3$ already in the same column.
 Additional comments to explain the translation:
 \begin{description}
\item [C:1] The two numbers are $(2\alpha_1)\alpha_0=6_4\alpha_0$ and $\alpha_0^2$.
\item [C:2] The subtractions will be
\begin{center}
\begin{tabular}{ccc}
$10^2$            &  10           &1\\ \hline
$3_2 $            & \bf 3       & \bf 4\\
                  & $-6_4\alpha_0$& $-\alpha_0^2$\\
\end{tabular}
\end{center}
\end{description}

\subsection{
 De Practica Geometrie  $\sqrt{6142}$ }
 The emphasized text is based on the transcript of Boncompagni \cite[p.20]{Boncompagni1862} and the translation by Hughes \cite[p.41-42]{Hughes2008}. See Figure \ref{fig:sqrt6142} for the corresponding line numbers.

 \begin{quotation}
 {\em Again, if you wish to find the root of 6142,}
 \begin{enumerate}
 \item {\em find the root of 61. It is 7
with remainder 12. Put the 7 under the 4 {\em{\color{blue}[L:1]}}}
\item {\em  and the 12 over 61 {\em{\color{blue}[L:2]}}.}
\item {\em Double 7 to
get 14. Put 4 under the 7 and 1 after the 7 {\em{\color{blue}[L:3]}}.}
\item {\em
You know about the union of the remainder 12 with 42 to make 1242 {\em{\color{blue}[C:1]}}.
As a four digit number three subtractions
are required. The first of these is the first digit
by the last digit. The second is the first digit by 4 that is before the
1. The third is the square of the first digit itself {\em{\color{blue}[C:2]}}.
We must gradually subtract
the three products from the four given digits so that the final product
falls in the first place. Because the number of digits exceeds the number
of products in the first figure, one digit must be joined to the last digits of
1242 with what follows, namely 12.
From this 12 subtract the first product.
Then the second product falls under the second place, and the last under
the first {\em{\color{blue}[C:3]}}.
Whence put 8 before 7 {\em{\color{blue}[L:4]}}.}
\item {\em
Multiply 8 by the given 1 and subtract
[the product] from 12 for a remainder of 4. Join this with the 4 that follows
in the second place to get 44.
Subtract from this the product of 8 and 4
(the 4 being under the 7) to yield a remainder of 12.
Place this over the 4{\color{red}2} {\em{\color{blue}[L:5]}}}
\item {\em Join it with 2 in the first place to make 122. Subtract from this 64 the
square of 8 for a remainder of 58 {\em{\color{blue}[L:6]}}.}
\end{enumerate}
{\em And thus you have 78 as the root of 6142
with remainder 58.}
 \end{quotation}

 \begin{figure}[ht!]
 \begin{minipage}[t]{.3\linewidth}
 \begin{center}
\fbox{\begin{tabular}{cccl}
        &       &       &$\left(58_6\right.$   \\
  $1_2$ &$2_2$  &$1_5$  & $2_5$  \\
 {\bf 6}&{\bf 1}&{\bf 4}&{\bf 2} \\
        &       &$7_1$  & $8_4$  \\
        &$1_3$  &$4_3$  &        \\
\end{tabular}}
\end{center}
\end{minipage}
\begin{minipage}[t]{.5\linewidth}
\(\begin{array}{l}
1:\ \alpha_1^2 \leq  \beta_3 10+\beta_2= \bm{61},  \alpha_1=7\\
2:\ \beta_3 10+\beta_2-\alpha_1^2=\bm{61}-49=12\text{ (in hundreds)}\\
3:\ 2\alpha_1=14 \text{ (in tens)}\\
4:\ 1_22_2\bm{42}-100\cdot\alpha_0-10\cdot 4\alpha_0-\alpha_0^2\geq 0, \alpha_0=8\\
5:\ ((12-8\cdot1)10+(4-4\cdot8)=12 \text{ (in tens)}\\
6:\ 1_52_5\bm2-8^2=58\\
\end{array}\)
\end{minipage}
 \caption{Fibonacci's example for  $\sqrt{6142}$ from De Practica Geometrie.}
 \label{fig:sqrt6142}
 \end{figure}
Additional comments.
\begin{description}
\item [C:1] Using the notation in the introductory section, $n=61$, $a=7=\alpha_1$, and $1_22_2=n-a^2$ so
\( 1_22_2 \cdot100+\beta_110+\beta_0= 1_22_2\bm{42} = N-(10\alpha_1)^2.\)
\item [C:2] The three numbers are
\(1_3\alpha_0, 4_3\alpha_0, \alpha_0^2.\) This will be \( 2\alpha_1\alpha_0\) (two numbers) and  \(\alpha_0^2\)
\item [C:3] To compute $2\alpha_1\alpha_0$ ($\alpha_0$ to be determined) consider the table
\begin{center}
\begin{tabular}{cccc}
$10^3$             &$10^2$            &  10           &1\\ \hline
$ 1_2$             &$2_2 $            &   \bf 4       & \bf 2\\
                   &$-1_3\alpha_0$    & $-4_3\alpha_0$& $-\alpha_0^2$\\
\end{tabular}
\end{center}
This will be \(1242-100\cdot\alpha_0-10\cdot 4\alpha_0-\alpha_0^2\)
which is supposed to be nonnegative. This is inequality (\ref{eq:bound_1}) and largest $\alpha_0=8$.
\item [C:4] For $\alpha_0=8$ we have
\begin{center}
\begin{tabular}{ccc}
$10^2$            &  10           &1\\ \hline
$ 1_22_2 $            &   \bf 4       & \bf 2\\
$-1_3\cdot 8_4$    & $-4_3\cdot 8_4$& $-8_4^2$\\\hline
$12-8_4\cdot 1_3=4$ \\
                   &$4\bm4 - 4_3\cdot8_4=12$\\
                   &                 &$12\bm2-8_4^2$
\end{tabular}
\end{center}
which is \((12-8\cdot1)10^2+(4-4\cdot8)10=120\) or 12 (in tens) and
\((12-8\cdot1)10^2+(4-4\cdot8)10+ 2-8^2 = 58.\)
\end{description}

 \subsection{
 De Practica Geometrie $\sqrt{8172}$ }
 The indented text is based on the transcript of Boncompagni \cite[p.20-21]{Boncompagni1862} and the translation by Hughes \cite[p.42]{Hughes2008}. See Figure \ref{fig:sqrt8172} for the corresponding lines.
 \begin{quotation}
 {\em If you wish to find the root of 8172}
 \begin{enumerate}
 \item {\em  put the root of 81 or 9 under the 7 {\em{\color{blue}[L:1]}}.}
 \item {\em  Double the 9 and put the 8 under the 9 and the 1 after it to the left} {\color{blue}[L:2]}.
 \item {\em Now the 1 and 8 must be multiplied by the first digit, one at a time. Then square the
first digit. And thus there are three products to be subtracted
gradually from 72, the remainder from the 81 after finding of the root of 81.
Whence, as we obviously know, nothing comes after it except 0 {\em{\color{blue}[L:3]}}.}
\item {\em  Since a step
is lacking, it is the first product that can be subtracted. Because if the first
product is subtracted from 7, the second needs be subtracted from 2.
But then there is no place from where to subtract the third product. Or in
another way: because the first place is a factor with any step, that step arises
from the multiplication. Since the product of the digit in the first place and
the digit in the third place, namely by 1, fills the third place, there is no place
for 72. Therefore the root of 8172 is 90 and the remainder is 72 {\em{\color{blue}[L:4]}}.}
\end{enumerate}
  \end{quotation}
  \begin{figure}[ht!]
 \begin{minipage}[t]{.3\linewidth}
   \begin{center}
\fbox{\begin{tabular}{cccl}
        &       &       &$\left(72_4\right.$   \\
 {\bf 8}&{\bf 1}&{\bf 7}&{\bf 2} \\
        &       &$9_1$  & $0_3$  \\
        &$1_2$  & $8_2$ &        \\
\end{tabular}}
\end{center}
\end{minipage}
\begin{minipage}[t]{.5\linewidth}
\(\begin{array}{l}
1:\ \alpha_1^2 \leq  \beta_3 10+\beta_2= \bm 8\bm 1,  \alpha_1=9\\
2:\ 2\alpha_1 =18 \text{ (in tens)}\\
3:\ \beta_1\cdot10 + \beta_0=\bm 7\bm2\\
3:\ 72-10(2\alpha_1)\alpha_0 - \alpha_0^2\geq 0, \alpha_0=0\\
4:\ 72-10(2 \alpha_1)\alpha_0 - \alpha_0^2 = 72
\end{array}\)
\end{minipage}
 \caption{Fibonacci's example for  $\sqrt{8172}$ from De Practica Geometrie}
 \label{fig:sqrt8172}
 \end{figure}

\subsection{
 De Practica Geometrie  $\sqrt{12345}$ }
 The indented text is based on the transcript of Boncompagni \cite[p.21]{Boncompagni1862} and the translation by Hughes \cite[p.43]{Hughes2008}. See Figure \ref{fig:sqrt12345} for the corresponding lines.
 \begin{quotation}{\em
For example: we wish to find the root of 12345. First, find the
root of 123 which is 11 with remainder 2. }
\begin{enumerate}
\item {\em
 Put the first 1 of 11 under the 3 and
the other 1 under the 4 {\em{\color{blue}[L:1]}}. }
\item { \em Double 11 to get 22 and place it under the 11 {\em{\color{blue}[L:2]}}.}
\item {\em Place that
remaining 2 above the 3 {\em{\color{blue}[L:3]}}. }
\item {\em Join it to the following digits to make 245, as a three digit
number. There will be the usual three multiplications. There will be a single
multiplication from each single digit. Hence, there needs be placed such a digit
before the 11 already in position. This will be multiplied by the first binomial,
then by the second, and finally by itself. The first multiplication can be taken
from the 2. Put what remains over the 3. And then another from the 4 and finally
another from the 5 which is in the first position {\em{\color{blue}[C:1]}}. And that digit is 1 {\em{\color{blue}[L:4]}}.}
\item {\em Having subtracted
the product of 1 by the first binomium and the product by the 2 over the
3, nothing remains. Likewise, having multiplied 1 by the following binomium and
subtracted it from 4, what remains is 2 over the 4. Having joined the 2 with the
5 in the first place, 25 is made. Having squared 1 and subtracted it from 25, 24
remains {\em{\color{blue}[C:2]} {\color{blue}[L:5]}}.}
\end{enumerate}
{\em  And thus you have 111 as the root of 12345 with remainder 24.}
\end{quotation}
\begin{figure}[th!]
 \begin{minipage}[t]{.4\linewidth}
 \begin{center}
\fbox{\begin{tabular}{ccccl}
        &        &       &      &$\left(24_5\right.$  \\
        &        &$2_3$  &      &        \\
 {\bf 1}&{\bf 2} &{\bf 3}&{\bf 4}&{\bf 5}\\
        &        &$1_1$ & $1_1$ & $1_4$  \\
        &        &$2_2$ & $2_2$ &        \\
\end{tabular}}
\end{center}
\end{minipage}
\begin{minipage}[t]{.55\linewidth}
\(\begin{array}{l}
1:\ \lfloor\sqrt{123}\rfloor = 11,  \text{ and insert the number 11}\\
1:\ \alpha_2=1, \alpha_1=1 \\
2:\ 2(\alpha_2 10 + \alpha_1)=22\text{ (in tens)}, \\
3:\ \beta_4 10^2+\beta_3 10 + \beta_2-(\alpha_2 10 + \alpha_1)^2 = \\
3:\  = 123-11^2=2  \text{ (in hundreds)}\\

4:\ 2_3{\bm 4 \bm 5}-22\cdot 10\alpha_0- \alpha_0^2\geq 0, \alpha_0=1\\
5:\ 2_3{\bm 4 \bm 5}-22\cdot 10- 1=24\\

\end{array}\)
\end{minipage}
 \caption{Fibonacci's example for  $\sqrt{12345}$ from De Practica Geometrie}
 \label{fig:sqrt12345}
 \end{figure}
 Additinal comments on the computation of the integer part of $\sqrt{12345}$.
 \begin{description}
 \item [C:1] The unknown (last digit) is $\alpha_0$ and we have
\begin{center}
\begin{tabular}{ccc}
$10^2$            &  10           &1\\ \hline
$ 2_3$            &   \bf 4       & \bf 5\\
$-2_2 \alpha_0$    & $-2_2\alpha_0$& $-\alpha_0^2$\\
\end{tabular}
\end{center}
or \(2_3{\bm 4 \bm 5}-22\cdot 10\alpha_0- \alpha_0^2\geq 0\) and largest $\alpha_0=1$.
\item [C:2] For  $\alpha_0=1$ we have
\begin{center}
\begin{tabular}{ccc}
$10^2$            &  10           &1\\ \hline
$ 2_3$            &   \bf 4       & \bf 5\\
$-2_2 \cdot 1$    & $-2_2\alpha_0$& $-\alpha_0^2$\\
= 0               & =2            & 4
\end{tabular}
\end{center}
which gives the remainder 24.
\end{description}
   \subsection{
De Practica Geometrie  $\sqrt{98765}$ }
The following is based on the transcript of Boncompagni \cite[p.21-22]{Boncompagni1862} and the translation by Hughes \cite[p.44]{Hughes2008}.
See Figure \ref{fig:sqrt98765}.
 \begin{quotation}{\em
 Again, if you wish to find the root of 98765, find first the root of 987. It
is 31 with a remainder of 26.}
\begin{enumerate}
\item {\em Put the 3 under the 7, the 1 under the 6} {\color{blue}[Line marked  1 in Figure \ref{fig:sqrt98765}]}.
\item {\em Put the remainder 26 over the 87, the 2 above the 8 and the 6 above the 7}
 {\color{blue}[L:2]}.
\item {\em Join 26 with the other two digits to make 2665 }{\color{blue}[C:1]}.
{\em Double the 31 and put the 6 under
the 3 and the 2 under the 1.} {\color{blue}[L:3]}.
\item  {\em And because four digits remain
in the number, we must gradually delete them by three multiplications that
have to be done with the first digit before 31. The multiplications are first by
6, second by 2, and the third is the square of the first digit} {\color{blue}[C:2]}.
{\em We must therefore
take the first product from 26. Whence we divide 26 by 6 (twice the 3
found in the partial root and placed under the 3) to get 4} {\color{blue}[C:3]}.
{\em Therefore the first
digit is 4 to be placed before the 31. So, put 4 before 31} {\color{blue}[L:4]}.
\item
{\em Since it can be done,
multiply 4 by 6. Subtract the product from 26 for a remainder of 2 that you
put over the 6 above the 7} {\color{blue}[L:5]}.
\item  {\em Join it with the 6 in the second place to make 26.
Then subtract from this number the product of 4 by 2  which placed
under the 1 leaves 18, 1 over the 2 and 8 over the 6} {\color{blue}[L:6] }.
\item {\em Join the 18 with the 5
in the first place to form 185. From this subtract the square of 4 to leave 169}
{\color{blue}[L:7]}.
\end{enumerate}
{\em Thus you have 314 for the root of 98765 with a remainder of 169} {\color{blue}[C:4]}.

    \end{quotation}
\begin{figure}[th!]
 \begin{minipage}[t]{.4\linewidth}
 \begin{center}
   \fbox{\begin{tabular}{ccccl}
        &        &       &      &$\left(169_7\right.$  \\
        &        &$1_6$  &      &        \\
        &        &$2_5$  &      &        \\
        & $2_2$  &$6_2$  &$8_6$ &        \\
 {\bf 9}&{\bf 8} &{\bf 7}&{\bf 6}&{\bf 5}\\
        &        &$3_1$ & $1_1$ & $4_4$  \\
        &        &$6_3$ & $2_3$ &        \\
  \end{tabular}}
 \end{center}
 \end{minipage}
\begin{minipage}[t]{.55\linewidth}
\(\begin{array}{l}
1\colon\ \lfloor\sqrt{987}\rfloor = 31,  \text{ and insert the number 31}\\
1:\ \alpha_2=3, \alpha_1=1 \\
2:\ \beta_4 10^2+\beta_3 10 + \beta_2-(\alpha_2 10 + \alpha_1)^2 = \\

2:\  =987-31^2=26,  \text{ (in hundreds)}\\
3:\ 2(\alpha_2 10 + \alpha_1)= 62\text{ (in tens)}, \\
4:\ \lfloor\frac{2_26_2}{6_3}\rfloor=4=\alpha_0\\
5:\ 2_2\cdot10 + (6_2-6\cdot4_4) = 26-24=2\text{ (in hundreds)}\\
6:\ 2_5 \cdot10 + \bm 6 - 2_3\cdot4_4  = 26-8=18 \text{ (in tens)}\\
7:\ 18_6\cdot 10 + \bm 5- 4_4^2=169
\end{array}\)
\end{minipage}
 \caption{Fibonacci's example for  $\sqrt{98765}$ from De Practica Geometrie}
 \label{fig:sqrt98765}
 \end{figure}
 Comments on computing the square root of 98765.
\begin{description}
 \item [C:1] The number \(2665=2_26_2{\bf 6}{\bf 5}\)
  \item[C:2] The four columns represents $10^3, 10^2, 10$ and $1$.
 \begin{center}
\begin{tabular}{cccc}
$10^3$ &$10^2$&  10      &     1    \\ \hline
$2_2$  &$6_2$ & ${\bm 6}$& ${\bm 5}$\\
       &$- 6_3\alpha_0 $&$- 2_3\alpha_0 $&$- \alpha_0^2$ \\

\end{tabular}
\end{center}
and the number is
\(2\cdot10^3 + (6-6\alpha_0)10^2 +(6-2\alpha_0)10 +5-\alpha_0^2.\)
 \item[C:3]Largest $\alpha_0$ in the inequality \(2665-62\alpha_0\cdot 10 - \alpha_0^2 \geq 0\) is $\alpha_0=4$. It follows from the inequality that
     \(\alpha_0 \leq \frac{2665}{620+\alpha_0}\). Possible approximations are
     \(\frac{2665}{620}\), \(\frac{266}{62}\), or \(\frac{26}{6}\) where the last approximation is used.
    \item[C:4] From comment [C:2] and $\alpha_0=4$ we have
    \begin{center}
\begin{tabular}{cccc}
$10^3$ &$10^2$&  10      &     1    \\ \hline
$2_2$  &$6_2$ & ${\bm 6}$& ${\bm 5}$\\
       &$- 6_3\cdot4_4 $&$- 2_3\cdot4_4 $&$- 4_4^2$ \\
\end{tabular}
\end{center}
Lines 5, 6 and 7 in Figure \ref{fig:sqrt98765} are to compute \(2665-62\alpha_0\cdot 10 - \alpha_0^2= 2665-620\cdot 4_4-4_4^2= 169\)
 \end{description}

   \subsection{
 De Practica Geometrie  $\sqrt{123456}$ }
 The indented text is based on the transcript of Boncompagni \cite[p.22]{Boncompagni1862} and the translation by Hughes \cite[p.44]{Hughes2008}.
 See Figure \ref{fig:sqrt123456}.
 \begin{quotation}
 {\em If you wish to find the root of a number with six digits, first find the root
of the last four digits and join the remainder with the following two digits, and
continue as before. For example, if we want to find the root of 123456, first find
the root of 1234, which is 35 with remainder 9.}
\begin{enumerate}
\item {\em  Put 35 under 45  {\em{\color{blue}[L:1]}}.}
\item {\em Put the remainder 9 above the 4 {\em{\color{blue}[L:2]}}}
\item {\color{blue}[L:3]}.
\item {\em Join 9 to the 56 to make 956. Because there
are three digits here and there are three multiplications. The first must be by
7, the second by zero, and the third the square of the number itself {\em{\color{blue}[C:1]}}.
Then we
know that the first multiplication must come out of the 9. So divide 9 by 7 to
produce 1, the first digit before 35 {\em{\color{blue}[L:4]}}.}
\item {\em
Multiply 1 by 7 and subtract the product
from 9 to leave the remainder 2 above the 9 {\em{\color{blue}[C:2]} {\color{blue}[L:5]}}.}
\item {\em Join the 2 to the 5 to make 25.
Multiply 1 by 0, subtract the product from 25, and 25 remains. Join this to the
following digit to make 256. From this subtract the square of 1 to leave 255 {\em{\color{blue}[L:6]}}.}
\end{enumerate}
{\em Thus we have 351 for the root of 123456 with remainder 255.}
    \end{quotation}
         \begin{figure}[ht!]
 \begin{minipage}[t]{.4\linewidth}
 \begin{center}
\fbox{\begin{tabular}{cccccl}
        &        &       &      &         &$\left(255_6\right.$  \\
       &         &       &$2_5$ &         &\\
        &        &       &$9_2$ &         &\\
 {\bf 1}&{\bf 2} &{\bf 3}&{\bf 4}&{\bf 5} &{\bf 6}\\
        &        &       & $3_1$ & $5_1$  &$1_4$   \\
        &        &       &$7_3$ & $0_3$   &      \\
\end{tabular}}
\end{center}
\end{minipage}
\begin{minipage}[t]{.55\linewidth}
\(\begin{array}{l}
1:\ \lfloor\sqrt{\bm{1234}}\rfloor = 35,  \text{ and insert 35}\\
1:\ \alpha_2=3, \alpha_1=5 \\
2:\ \beta_5 10^3+\beta_4 10^2+\beta_3 10 + \beta_2-(\alpha_2 10 + \alpha_1)^2 =\\
2:\  =\bm{12345}-35^2=9  \text{ (in hundreds) }\\
3:\ 2(\alpha_2 10 + \alpha_1)=70 \text{ (in tens)}, \\
4:\ \alpha_0\approx\lfloor\frac{9_2}{7_3}\rfloor, \alpha_0=1\\
5:\ 9_2-7_3\cdot1_4=2 \text{ (in hundreds)}\\
6:\ 2_5\cdot 10 + \bm 5 - 0_3\cdot1_4 =25 \text{ (in tens)}\\
6:\ 25\cdot10+\bm 6 - 1_4^2=255
\end{array}\)
\end{minipage}
 \caption{Fibonacci's example for  $\sqrt{123456}$ from De Practica Geometrie}
 \label{fig:sqrt123456}
 \end{figure}
 Additional comments on the square root computation.
 \begin{description}
   \item[C:1] \

 \begin{center}
\begin{tabular}{ccc}
  $10^2$&  10      &     1    \\ \hline
  $9_1$ &   $\bm5$ & ${\bm 6}$\\
  $- 7_3\alpha_0 $&$- 0_3\alpha_0 $&$- \alpha_0^2$ \\
\end{tabular}
\end{center}
   \item[C:2] \
   \begin{center}
\begin{tabular}{ccc}
  $10^2$&  10      &     1    \\ \hline
  $9_1$ &   $\bm5$ & ${\bm 6}$\\
  $- 7_31_4$&$- 0_3\cdot1_4 $&$- 1_4^2$ \\
\end{tabular}
\end{center}
 \end{description}
   \subsection{
De Practica Geometrie $\sqrt{987654}$ }
 The indented text is based on the transcript of Boncompagni \cite[p.22]{Boncompagni1862} and the translation by Hughes \cite[p.45]{Hughes2008}. In the case of no subscripts on the inserted digits in the table in Figure \ref{fig:sqrt987654}, the table can be reduced when the same digit are consecutive in a column.
 \begin{quotation}
 {\em   Again, if you wish to find the root of 987654, find first the root of 9876 which
is 99 with a remainder of 75.}
\begin{enumerate}
\item {\em  Put 99 under 65 {\em{\color{blue}[L:1]}}.}
\item {\em  Put 75 above 76 {\em{\color{blue}[L:2]}}.}
\item {\em Double the root you found to get 198.
From this number, put 8 under the 9 in
the second place, 9 under the 9 in the third place, and 1 after the last {\em{\color{blue}[L:3]}}.}
\item {\em Join 75 with 54 to make 7554 {\em{\color{blue}[C:1]}}.
There are four digits in this number which we must delete
gradually by four multiplications by the three digits 1, 9, and 8 with the first digit
and finally by the square of itself. Whence the first product must be subtracted
only from 7.
If we know that the first digit is 3 {\em{\color{blue}[L:4]}{\color{blue}[C:2]}}.}
\item {\em
The product of 3 and 1 in 198
subtracted from 7 leaves 4 above 7 itself {\em{\color{blue}[L:5]}}.}
\item {\em Join it to 5 in the third place to make 45.
Now subtract from this the product of 3 and the 9 in 198, to have 18 remain over
the fourth and third places{\em{\color{blue}[L:6]}}.}
\item {\em
Join 18 to 5 in the second place to make 185. Subtract
from this the product of 3 and the 8 in 198. 161 remains above the fourth, third,
and second places {\em{\color{blue}[L:7]}}.}
\item {\em
Join 161 to 4 in the first place to make 1614. Subtract the
square of 3 to give a remainder of 1605 {\em{\color{blue}[L:8]}}.}
\end{enumerate}
{\em And so we have 993 as the root.}
    \end{quotation}
\begin{figure}[ht!]
 \begin{minipage}[t]{.4\linewidth}
 \begin{center}
\fbox{\begin{tabular}{cccccl}
        &        & $1_7$ &       &        &$\left(1605_8\right.$  \\
        &        & $1_6$ &$6_7$  &        &\\
        &        & $4_5$ &$8_6$  &        &\\
        &        & $7_2$ &$5_2$  &{\color{red}$1_7$}        &\\
 {\bf 9}&{\bf 8} &{\bf 7}&{\bf 6}&{\bf 5} &{\bf 4}\\
        &        &       & $9_1$ & $9_1$  &$3_4$   \\
        &        &$ 1_3$ & $9_3$ & $8_3$  &      \\
\end{tabular}}
\end{center}
\end{minipage}
\begin{minipage}[t]{.55\linewidth}
\(\begin{array}{l}
1:\ \lfloor\sqrt{9876}\rfloor = 99,  \text{ and insert the number 99}\\
1:\ \alpha_2=9, \alpha_1=9 \\
2:\ \beta_4 10^2+\beta_3 10 + \beta_2-(\alpha_2 10 + \alpha_1)^2 = \\
2:\ 9876-99^2=75  \text{ (in hundreds)}\\
3:\ 2(\alpha_2 10 + \alpha_1)=198\text{ (in tens)}, \\
4:\ \alpha_0=3\\
5:\ 7_2- 1_3\cdot 3_4=7-3=4 \text{ (in thousands)}\\
6:\ 4_55_2 - 9_3\cdot 3_4=45-27=18 \text{ (in hundreds)}\\
7:\ 1_68_6\bm5-8_3\cdot 3_4=185-24=161 \text{ (in tens)}\\
8:\ 1_76_71_7\bm4- 3_4^2=1605
\end{array}\)
\end{minipage}
 \caption{Fibonacci's example for  $\sqrt{987654}$ from De Practica Geometrie}
 \label{fig:sqrt987654}
 \end{figure}
 Additinal comments on computing the integer part of $\sqrt{987654}$.
 \begin{description}
   \item[C:1] \(7_25_2\bm5\bm4=7554\)
 \item[C:2] \(\alpha_0\approx\lfloor\frac{7_25_2\bm 5}{1_39_38_3}\rfloor=3,\) or
 \(\alpha_0\approx\lfloor\frac{7_25_2}{1_39_3}\rfloor=3 \) or largest $\alpha_0$ so that
 \[7554-198\cdot10\cdot \alpha_0- \alpha_0^2\geq 0.\]
 \item[C:3] Need to compute the remainder
   \begin{center}
\begin{tabular}{cccc}
$10^3$ &$10^2$&  10      &     1    \\ \hline
$7_2$  &$5_2$ & ${\bm 5}$& ${\bm 4}$\\
$- 1_3\alpha_0 $&$- 9_3\alpha_0 $&$- 8_3\alpha_0 $&$- \alpha_0^2$ \\
\end{tabular}
\end{center}
 For $\alpha_0=3$ the table is
 \begin{center}
\begin{tabular}{c|c|c|c}
$10^3$ &$10^2$&  10      &     1    \\ \hline
$7_2$  &$5_2$ & ${\bm 5}$& ${\bm 4}$\\
$- 1_3\cdot 3_4 $&$- 9_3\cdot 3_4 $&$- 8_3\cdot 3_4 $&$- 3_4^2$\\ \hline
\(7_2- 1_3\cdot 3_4\)&
\(4_55_2 - 9_3\cdot 3_4\)& \(1_68_6\bm5-8_3\cdot 3_4\) &$1_76_71_7\bm4- 3_4^2$ \\
$=7-3=4_5$ & $=45-27=1_68_6$&$=185-24=1_76_71_7$ &$=1605 $
\end{tabular}
\end{center}
   \end{description}

\subsection{
De Practica Geometrie  $\sqrt{9876543}$ }
The indented text is based on the transcript of Boncompagni \cite[p.22]{Boncompagni1862} and the translation by Hughes \cite[p.45--46]{Hughes2008}. If there are no subscripts on the digits in the table, the table can be reduced when the same digit are consecutive in a column (See Figure \ref{fig:sqrt9876543}).
 \begin{quotation}
 \emph{ Likewise, if you wish to find the root of a seven digit number, find first the
root of the last five digits. Join any remainder with the remaining two digits,
and proceed as before. For example: we want to find the root of 9876543. Since the root of the seven digit number has four digits, we must
put the first digit of the root under the fourth place, the following digit under
the third, and the next one under the second place. You will find these three
digits by finding the root of the last five digits. }
\begin{enumerate}
\item {\em
Their root is 314 }{\color{blue}[L:1]}.
\item {\em The remainder
169 you place above the fifth, fourth, and third places{\em{\color{blue}[L:2]}}}.
\item {\em
Join 169 to 43 the
remaining two digits to make 16943. Then double the 314 that you found and
put the 8 under the 4, the 2 under the 1, and the 6 under the 3. {\em \color{blue}[L:3]}}
\item {\em  You will multiply
the first digit by those three digits, and then square the first digit. Consequently
there are four products which we must subtract step-by-step from
the five digit figure 16943.
So the first product to be subtracted is from the last
two figures, namely from 16. Hence divide 16 by 6 the first product to
get 2 {\em{\color{blue}[L:4]}}.}
\item{\em
Put the 2 in the first place of the root. {\em \color{blue}[C:1]}
Then multiply it by the 6 of 628
to get 12 that you subtract from 16 for a remainder of 4. Put this over the 6 {\em \color{blue}[L:5]}.}
\item{\em Join 4 with 9 to make 49. From this subtract the product of 2 by 2 to leave 45
over the 49 {\em \color{blue}[L:6]}.}
\item {\em
Join this with the 4 in the second place to make 454. Subtract from
this the product of 2 by 8 to leave 438.
Put the 4 over the fourth place, 3 above the
third, and 8 above the second {\em \color{blue}[L:7]}.}
\item{\em
Join 438 with the 3 in the first place to make
4383. From this subtract the square of 2 for a remainder of 4379 {\em \color{blue}[L:8]}.
This is less
than twice the root you found. Consequently, the root is 3142 {\em \color{blue}[C:2]}.}
\end{enumerate}
\end{quotation}
     \begin{figure}[ht!]
 \begin{minipage}[t]{.4\linewidth}
 \begin{center}
\fbox{\begin{tabular}{ccccccl} &        &       &      &         &       &$\left(4379_8\right.$  \\
        &        &       &{\color{red} $4_7$} &   &    &  \\
        &        &       &{\color{red} $4_6$} &{\color{red} $3_7$}    &    &  \\
        &        &       &{\color{red} $4_5$} &{\color{red} $5_6$}   &    &  \\
        &        & $1_2$ &$6_2$               &  $9_2$               &{\color{red} $8_7$}        &\\
 {\bf 9}&{\bf 8} &{\bf 7}&{\bf 6}&{\bf 5} &{\bf 4}&{\bf 3}\\
        &        &       & $3_1$ & $1_1$  &$4_1$    &$2_4$\\
        &        &       &{\color{red} $6_3$} &{\color{red} $2_3$}  &{\color{red}$8_3$}    &\\
\end{tabular}}
\end{center}
\end{minipage}
\begin{minipage}[t]{.55\linewidth}
\(\begin{array}{l}
1:\ 1:\ \lfloor\sqrt{98765}\rfloor =314  \text{ and insert the number 314}\\
1:\ \alpha_3=3, \alpha_2=1, \alpha_1=4 \\
2:\ \beta_6 10^4+\beta_5 10^3+\beta_4 10^2 + \beta_310-\beta_2 \\
2:\ \quad -(\alpha_3 10^2+\alpha_2 10 + \alpha_1)^2 =\\
2:\ \quad 98765-314^2=169  \text{ (in hundreds)}\\
3:\ 2(\alpha_3 10^2+\alpha_2 10 + \alpha_1)=2\cdot314=628\\
4:\ \alpha_0\approx\lfloor\frac{1_26_2}{6_3}\rfloor, \alpha_0=2\\
5:\ 1_26_2-6_3\cdot2_4 =16-12=4 \text{ (in thousands)}\\
6:\ 4_5 9_2 - 2_3\cdot2_4 =49-4=45 \text{ (in hundreds)}\\
7:\ 4_6 5_6 \bm 4 - 8_3\cdot2_4 =454-16=438 \text{ (in tens)}\\
8:\ 4_7 3_7 8_7 \bm 3 - 2_4^2 =4383-16=4379 \text{ (in ones)}
\end{array}\)
\end{minipage}
 \caption{Fibonacci's example for  $\sqrt{9876543}$ from De Practica Geometrie}
 \label{fig:sqrt9876543}
 \end{figure}
 Additional comments:
 \begin{description}
 \item [C:1]\ \begin{center}
\begin{tabular}{c|c|c|c}
$10^3$ &$10^2$&  10      &     1    \\ \hline
$1_26_2$ & $9_2$ &${\bm 4}$& ${\bm 3}$\\
$- 6_3\cdot2_4 $&$- 2_3\cdot2_4 $&$- 8_3\cdot2_4 $&$- 2_4^2$ \\ \hline
$1_26_2- 6_3\cdot2_4 =$&$4_5 9_2 - 2_3\cdot2_4 =$& \(4_6 5_6 \bm 4 - 8_3\cdot2_4 =\) &
\(4_7 3_7 8_7 \bm 3 - 2_4^2 =\)\\
$=16-12=4_5$&$=49-4=4_65_6$&$ =454-16=4_73_78_7$&$=4383-16=4379$
\end{tabular}
\end{center}
 \item[C:2] This is inequality (\ref{eq:bound_2})
 \(4379 \leq 2\cdot3142\)
Let $a$ be the integer part of the square root of the integer $N$, then  $N-a^2\geq 0$ and $N-(a+1)^2= N-a^2 -2a-1<0$ which is inequality (\ref{eq:bound_2}).
\end{description}

\section{Computing the fractional part}\label{sec:4}

Let $a$ be an approximation to the square root of positive integer $N$ (not a perfect square) and let $r$ be the residual $r=N-a^2$. A common approximation is
\begin{equation}\label{eq:iter_1}
\sqrt{N} = \sqrt{a^2+r},\  a_1=a+\frac{r}{2a}.
\end{equation}
The residual (or remainder) is a measure on how good the solution is. If the residual $|r|$ is not sufficiently small, the iteration is repeated
\begin{equation}\label{eq:iter_2}
a_2 = a_1+ \frac{r_1}{2a_1} \text{ where } r_1= N-a_1^2.
\end{equation}
Further, from the equation defining $a_1$
\[a_1^2-N= \left(\frac{N-a^2}{2a}\right)^2, \]
so the new residual is the square of the correction \(\frac{r}{2a}\) of $a$. This makes the hand computation of the fractional part very efficient.
\begin{equation}\label{eq:iter}
a_1=a+\frac{r}{2a},\  a_2=a_1-\frac{(r/2a)^2}{2a_1}.
\end{equation}
The residual $r$ is the sum of the two last rows in the tables in Liber Abaci and is on the last row in the tables in
De Practica Geometrie (except the last digit) so the computational procedure used by Fibonacci is both very 
natural and efficient.

In a digit-by-digit method, common choices are $a=\lfloor\sqrt{N}\rfloor$ or $a=\lfloor\sqrt{N}\rfloor+1$
In the first case the residual $r$ is positive, while in the other case the residual will be negative.  All further generated residuals  $r_i,\ i\geq 1$ are negative. In Liber Abaci Fibonacci chooses $a=\lfloor\sqrt{N}\rfloor$ in the five examples discussed in this paper.

The same iterates are generated by Heron's method (also called the Babylonial method) and the Newton--Raphson method. Heron of Alexandra (active around 60 AD) computed the iterates as \cite[p.323-324]{Heath1921}
\begin{equation}\label{eq:iter_Heron}
a_1=\frac{1}{2}\left(a+\frac{N}{a}\right), \ a_2=\frac{1}{2}\left(a_1+\frac{N}{a_1}\right).
\end{equation}
Newton-Raphson's method applied to the equation $p(x)=x^2-N$ is
\begin{equation}\label{eq:iter_Newton}
 a_{i+1}= a_i-\frac{p(a_i)}{p'(a_i)}= a_i+ \frac{N-a_i^2}{2a_i}
 \end{equation}
and for $a_0=a$ the three methods give the same iterates but by computing by hand the methods are different.

In his work on approximation methods of Fibonacci, Glushkov 1976\cite{Glushkov1976} points out that (\ref{eq:iter}) is used without giving any examples.

Cantor in 1900 \cite{Cantor1900} points out for  $a^2<N<(a+1)^2$ then the two approximations are
\[\sqrt{N} \approx a+\frac{N-a^2}{2a}, \text{ and }
  \sqrt{N} \approx a+\frac{N-a^2}{2a} - \frac{\left(\frac{N-a^2}{2a}\right)^2}{2(a+\frac{N-a^2}{2a})}\]
As an example Cantor takes the example $\sqrt{927435}$ where $a=\lfloor\sqrt{927435}\rfloor=963$ and $N-a^2=66$.
\[\sqrt{927435}\approx 963+\frac{11}{321}- \frac{\left(\frac{11}{321}\right)^2}{2(963+\frac{11}{321})}.\]

The first example Fibonacci gives in the text is to demonstrate how to compute the fractional part of $\sqrt{10}$
\cite[490--491]{Sigler2002} \cite[p.353]{Boncompagni1857}\cite[p.548]{Giusti2020}.
\begin{quote}{\em You find the largest (square) root in 10 that can
be found which is the root of 9, and the 9 you subtract from the 10 leaving 1;
you divide the 1 by double the found root, namely 6; the quotient is $\frac{1}{6}$ which you add to the found 3; there will be $\frac{1}{6}3$ that is slightly larger than the root of
10, because when $\frac{1}{6}3$ is multiplied by itself this makes $\frac{1}{36}$ plus the number ten.}
\end{quote}
For \(a=\lfloor \sqrt{N}\rfloor\) where $N=10$, the integer part of $\sqrt{N}$ is $a=3$. The correction is \(\frac{N-a^2}{2a}=\frac{1}{6}\) so the first approximation is \(a+ \frac{N-a^2}{2a}=3\frac{1}{6}\) and $a^2= 10\frac{1}{36}$. For the next iteration
(\cite[p.491]{Sigler2002}, \cite[p.353]{Boncompagni1857}\cite[p.548]{Giusti2020}):
\begin{quote}{\em Suppose you wish to find something closer to the root of the 10, then you
divide the $\frac{1}{36}$ by double the $\frac{1}{6}3$ yielding $\frac{1}{228}$ which you subtract from the $\frac{1}{6}3$,
and you will have the proposed closer number.}
\end{quote}
For $a_1= 3\frac{1}{6}$, so the new approximation is \(a_2=a_1 - \frac{(\frac{1}{6})^2}{2a_1}= 3\frac{1}{6} - \frac{1}{228}\).

The following list follows the order of derivation of the integer part of the square root in Section \ref{sec:LiberAbaci}.
\begin{itemize}
\item Section 3.1 demonstrates computing $\lfloor\sqrt{743}\rfloor$, first determine the integer part 
$a=\lfloor\sqrt{743}\rfloor=27$ with residual $N-a^2=14$.
    \begin{quote}{\em \ldots  the 14 you divide by double the 27, or you divide half of 14, namely 7, by the 27;
    the quotient will be $\frac{7}{27}$ which you add to the found 27; there will be
$27\frac{7}{27}$ for the root of 743}
\end{quote}
Let $N=743, a=\lfloor\sqrt{743}\rfloor=27$ and $N-a^2=14$. The approximation is \(a+ \frac{N-a^2}{2a}=27+\frac{14}{2\cdot27}=27\frac{7}{27}\).
This approximation is in \cite[p.158]{Friedlein1869} and Pisano and Bussotti \cite[p.143--144]{Pisano2015}.
\item Section 3.2 demonstrates computing $\lfloor\sqrt{8754}\rfloor$ and the residual.
\begin{quote}{\em \ldots therefore the root of 8754 is in integers 93, and 105 remains;
this you divide by the double of the 93; the quotient will be $\frac{35}{62}$ which you add to the found 93; there will be $93\frac{35}{62}$  for the root of 8754.}
\end{quote}
Let $N=8754,  a=\lfloor\sqrt{N}\rfloor=93$ and $N-a^2= 105$. The approximation is
    \(a+ \frac{N-a^2}{2a}=93+\frac{105}{2\cdot93}=93\frac{35}{62}\).
 This approximation is in \cite[p.158]{Friedlein1869}.
 
\item Section 3.3 demonstrates computing $\lfloor\sqrt{12345}\rfloor$ and the residual.
\begin{quote}
{\em \ldots thus you will have a number of three figures, namely 111, for the root of 12345, as should be, and there remains the 24 beyond the root; half of
this, namely 12, you divide by the 111; the quotient will be $\frac{4}{37}$ which added to
the 111 yields $111\frac{4}{37}$ for the root of 12345.}
\end{quote}
$N=12345,  a=\lfloor\sqrt{N}\rfloor=111$ and $N-a^2= 24$. The approximation is
    \(a+ \frac{N-a^2}{2a}=111+\frac{24}{2\cdot111}=111\frac{4}{37}\).
    This approximation is in \cite[p.158]{Friedlein1869}.
    
\item Section 3.4 demonstrates computing $\lfloor\sqrt{927435}\rfloor= 963$ and the residual 66.
\begin{quote}\ldots{\em  there remains 66 of which half, namely 33, you divide by the 963; the quotient will be $\frac{11}{321}$,
and thus you will have $963\frac{11}{321}$ for the sought root, and
the product of it by itself yields a result greater than the sought number by the
amount of the multiplication of the fraction by itself, namely $\frac{11}{321}$.}
\end{quote}
$N=927435,  a=\lfloor\sqrt{N}\rfloor=963$ and $N-a^2= 66$. The approximation is
    \(a+ \frac{N-a^2}{2a}=963+\frac{66}{2\cdot963}=963\frac{11}{321}\).
    \begin{quote}
    {\em Therefore
if you wish to find a number closer to the root of 927435, then you multiply the $\frac{11}{321}$ by itself, and that which will result you divide by double the $963\frac{11}{321}$; the quotient will be less than the abovewritten.}
    \end{quote}
    $N=927435, a= 963\frac{11}{321}$, $a^2- N= \left(\frac{11}{321}\right)^2$ and subtract the correction $\frac{N-a^2}{2a}$.
   Cantor \cite[p.30]{Cantor1900} and Friedlein \cite[p.158]{Friedlein1869} use the example $\sqrt{927435}$
\[\sqrt{927435}\approx 963\frac{11}{321}- \frac{\left(\frac{11}{321}\right)^2}{2(963\frac{11}{321})}.\]

\item  Section 3.5 demonstrates computing $\lfloor\sqrt{72340000}\rfloor=8505$ and the residual 4975.  
    \begin{quote}
    {\em \ldots in the generated order you
divide the 4975 by the 17010 yielding about one fourth; and thus for the root
of 72340000 you have $\frac{1}{4}8505$ which you divide by the 100; the quotient will be
$\frac{1}{400}\frac{1}{20}85$ for the root of 7234.}
    \end{quote}
    $N=72340000,  a=\lfloor\sqrt{N}\rfloor=8505$ and $N-a^2= 4975$. Thus
    $\frac{N-a^2}{2a} \approx \frac{1}{4}$. The approximation of $\sqrt{72340000}$ is then $8505\frac{1}{4}$ and the approximation of
    $\sqrt{7234}$ is $85 + \frac{1}{20}+ \frac{1}{400}$
\end{itemize}
Liber Abaci contains more square root computation and computing the fractional part, but details are missing.
As an example we find $\sqrt{7}$ a little less than $2\frac{2}{3}$  which can only be explained using $a=\lfloor\sqrt{N}\rfloor+1$.
This and other examples in Liber Abaci and De Practica Geometrie are discussed in \cite{Steihaug2024}. There is no evidence that Fibonacci has used other than  two formulas (\ref{eq:iter}) and a rounding for computing 
an approximation of the fractional part of the square root and in all the ``worked out" examples in Liber Abaci 
the initial approximation is $a=\lfloor\sqrt{N}\rfloor$. 
\section{Concluding remarks}
In this paper we have shown that the verbal description of computing the integer part of the square root of a natural number in Fibonacci's Liber Abaci and De Practica Geometrie is uniquely described by a table and corresponding calculation. Computing the integer part of the square root and the fractional part in Liber Abaci 
takes the shape of a textbook in mathematics. The MSS, transcriptions and translations contain minor misalignment or missing numbers that are easily corrected. Further we show that Fibonacci is consistent in Liber Abaci in computing the fractional part and
the method is computational different from Heron's method to compute an approximation.

\bibliographystyle{plain}
\bibliography{Fibonacci,../Test_bibliography_test}
\end{document}